\documentclass[final]{siamltex}
\usepackage{pstricks}
\usepackage{amssymb}
\usepackage{amsmath}
\usepackage{amscd}
\usepackage{latexsym}
\usepackage{cite}
\usepackage{epic}
\usepackage{eepic}
\usepackage{fancybox}
\usepackage{longtable}
\usepackage{nicefrac}
\usepackage[dvips]{graphicx}
\usepackage{psfrag}
\usepackage{showkeys}
\usepackage{ifthen}

\usepackage[dvips,bookmarks=true]{hyperref}

\title{Self-adjoint $\curl$ Operators}
\author{Ralf Hiptmair\thanks{Seminar for Applied Mathematics, ETH Zurich, CH-8092
    Zurich, hiptmair\symbol{64}sam.math.ethz.ch}
  \and  Peter Robert Kotiuga\thanks{Boston University, Dept. of EC. Eng., prk@bu.edu}
  \and S{\'e}bastien Tordeux\thanks{Institut de Math{\'e}matiques de Toulouse, INSA-Toulouse, sebastien.tordeux@insa-toulouse.fr}}
\date{\today}

\providecommand{\Bdsymb}{\partial}                   
\providecommand{\Div}{\operatorname{div}}          
\providecommand{\bDiv}{{\Div}_{\Bdsymb}}           
\providecommand{\curl}{\operatorname{{\bf curl}}}  
\providecommand{\bcurl}{{\curl}_{\Bdsymb}}         
\providecommand{\bscurl}{\operatorname{curl}_{\Bdsymb}} 
\renewcommand{\grad}{\operatorname{{\bf grad}}}       
\providecommand{\bgrad}{\grad_{\Bdsymb}}                

\providecommand{\Dim}{\operatorname{dim}}            
\providecommand{\dim}{\Dim}

\providecommand*{\Kern}[1]{\mathcal{N}({#1})} 

\renewcommand*{\Range}[1]{\mathcal{R}({#1})} 

\newcommand{\Ve}{{\mathbf{e}}}

\newcommand{\Vn}{{\mathbf{n}}}

\newcommand{\Vu}{{\mathbf{u}}}
\newcommand{\Vv}{{\mathbf{v}}}

\newcommand{\Vx}{{\mathbf{x}}}

\providecommand{\Bx}{{\boldsymbol{x}}}
\providecommand{\By}{{\boldsymbol{y}}}

\newcommand{\VC}{{\mathbf{C}}}

\newcommand{\VG}{{\mathbf{G}}}
\newcommand{\VH}{{\mathbf{H}}}
\newcommand{\VI}{{\mathbf{I}}}
\newcommand{\VJ}{{\mathbf{J}}}

\providecommand{\BL}{{\boldsymbol{L}}}

\newcommand{\nablabf}{\boldsymbol{\nabla}}

\newcommand{\varphibf}{\boldsymbol{\varphi}}

\newcommand{\Phibf}{\boldsymbol{\Phi}}

\providecommand{\Cd}{{\cal D}}

\providecommand{\Ch}{{\cal H}}

\providecommand{\bbC}{\mathbb{C}}

\providecommand{\bbH}{\mathbb{H}}

\providecommand{\bbN}{\mathbb{N}}

\providecommand{\bbR}{\mathbb{R}}

\providecommand{\bbZ}{\mathbb{Z}}

\providecommand*{\N}[1]{\left\|{#1}\right\|} 
\newcommand*{\SN}[1]{\left|{#1}\right|}      

\newcommand{\defaultdomain}{D}

\providecommand*{\Lp}[2][\defaultdomain]{L^{#2}({#1})}
\newcommand*{\Lpv}[2][\defaultdomain]{\BL^{#2}({#1})}
\newcommand*{\NLp}[3][\defaultdomain]{\N{#2}_{\Lp[#1]{#3}}}

\newcommand*{\Ltwo}[1][\defaultdomain]{\Lp[#1]{2}}
\newcommand*{\Ltwov}[1][\defaultdomain]{\Lpv[#1]{2}}
\newcommand*{\NLtwo}[2][\defaultdomain]{\NLp[#1]{#2}{2}}

\newcommand*{\Hm}[2][\defaultdomain]{H^{#2}({#1})}

\newcommand*{\Hone}[1][\defaultdomain]{\Hm[#1]{1}}

\newcommand*{\NHone}[2][\defaultdomain]{{\N{#2}}_{\Hone[{#1}]}}

\newcommand*{\SNHone}[2][\defaultdomain]{{\SN{#2}}_{\Hone[{#1}]}}

\newcommand{\hlb}{\frac{1}{2}}
\providecommand*{\Hh}[1][\defaultboundary]{\Hm[#1]{\hlb}}
\newcommand*{\Hhv}[1][\defaultboundary]{\Hmv[#1]{\hlb}}

\newcommand{\vsobsymb}{\boldsymbol{H}}
\newcommand*{\Hdiv}[1][\defaultdomain]{\vsobsymb(\Div,{#1})}

\newcommand*{\Hcurl}[1][\defaultdomain]{\vsobsymb(\curl,{#1})}
\newcommand*{\bHcurl}[2][\defaultdomain]{\vsobsymb_{#2}(\curl,{#1})}
\newcommand*{\zbHcurl}[1][\defaultdomain]{\bHcurl[#1]{0}}

\newcommand*{\NHcurl}[2][\defaultdomain]{\N{#2}_{\Hcurl[#1]}}

\newcommand{\sH}{\mathsf{H}}

\newcommand{\sR}{\mathsf{R}}
\newcommand{\sS}{\mathsf{S}}
\newcommand{\sT}{\mathsf{T}}

\newcommand{\su}{\mathsf{u}}
\newcommand{\sv}{\mathsf{v}}

\newcommand{\sZ}{\mathsf{Z}}
\newcommand{\sD}{\mathsf{D}}

\newcommand*{\dom}{\defaultdomain}
\newcommand*{\zHcurl}[1][\dom]{\bHcurl[#1]{0}}
\newcommand*{\sd}{\operatorname{\textsf{d}}}
\newcommand*{\sId}{\mathsf{Id}}
\newcommand*{\hop}{\star}
\newcommand*{\symporth}{\sharp}
\newcommand*{\tscurl}[1][\partial\dom]{W^{-\frac{1}{2},1}(\sd,#1)}
\newcommand*{\tsgrad}[1][\partial\dom]{W^{-\frac{1}{2},0}(\sd,#1)}

\newcounter{remcnt}
\newenvironment{Remark}{\refstepcounter{remcnt}\textit{Remark \arabic{remcnt}}.}{}

\begin{document}

\maketitle

\centerline{\textsf{Report 2008-27, SAM, ETH Z\"urich, 
\url{http://www.sam.math.ethz.ch/reports/}}}

\begin{abstract}
  We study the exterior derivative as a symmetric unbounded operator on square
  integrable 1-forms on a 3D bounded domain $D$. We aim to identify boundary conditions
  that render this operator self-adjoint. By the symplectic version of the
  Glazman-Krein-Naimark theorem this amounts to identifying complete Lagrangian
  subspaces of the trace space of $\Hcurl$ equipped with a symplectic pairing arising
  from the $\wedge$-product of 1-forms on $\partial D$. Substantially generalizing earlier results,
  we characterize Lagrangian subspaces associated with closed and co-closed traces.
  In the case of non-trivial topology of the domain, different contributions from
  co-homology spaces also distinguish different self-adjoint extension. Finally, all
  self-adjoint extensions discussed in the paper are shown to possess a discrete
  point spectrum, and their relationship with $\curl\curl$-operators is discussed.
\end{abstract}



\pagestyle{myheadings}
\thispagestyle{plain}
\markboth{R. Hiptmair, P.R. Kotiuga, S. Tordeux}{Self-adjoint $\curl$ operators}


\section{Introduction}
\label{sec:introduction}

The $\curl$ operator is pervasive in field models, in particular in electromagnetics,
but hardly ever occurs in isolation. Most often we encounter a $\curl\curl$ operator
and its properties are starkly different from those of the $\curl$ alone. We devote
the final section of this article to investigation of their relationship.

The notable exception, starring a sovereign $\curl$, is the question of stable
force-free magnetic fields in plasma physics. They are solutions of the eigenvalue
problem
\begin{gather}
  \label{eq:63}
  \alpha\in\bbR\setminus\{0\}:\quad\curl\VH = \alpha\VH\;,
\end{gather}
posed on a suitable domain, see \cite{LUN50,CHK57,JET70,PIC96}. A solution theory for
\eqref{eq:63} must scrutinize the spectral properties of the $\curl$ operator. The
mature theory of unbounded operators in Hilbert spaces is a powerful tool. This
approach was pioneered by R. Picard \cite{PIC76,PIC98a,PIC98}, see also \cite{YOG90}.

The main thrust of research was to convert $\curl$ into a self-adjoint operator
by a suitable choice of domains of definition. This is suggested by the following
Green's formula for the $\curl$ operator:
\begin{gather}
  \label{eq:curlGreen}
  \int\nolimits_{\dom} \curl\Vu\cdot\Vv - \curl\Vv\cdot\Vu\,\mathrm{d}\Vx =
  \int\nolimits_{\partial\dom} (\Vu\times\Vv)\cdot\Vn\,\mathrm{d}S\;,
\end{gather}
for any domain $\dom\subset\bbR^{3}$ with sufficiently regular boundary
$\partial\dom$ and $\Vu,\Vv\in C^{1}(\overline{\dom})$. This reveals that
the $\curl$ operator is truly symmetric, for instance, when acting on
vector fields with vanishing tangential components on $\partial \dom$.

On bounded domains $\dom$ several instances of what qualifies as a self-adjoint
$\curl$ operators were found. Invariably, their domains were defined through
judiciously chosen boundary conditions. It also became clear that the topological
properties of $\dom$ have to be taken into account carefully, see
\cite[Thm.~2.4]{PIC98a} and \cite[Sect.~4]{YOG90}.

In this paper we carry these developments further with quite a few novel twists: we
try to give a rather systematic treatment of different options to obtain self-adjoint
$\curl$ operators. It is known that the $\curl$ operator is an incarnation of the
exterior derivative of 1-forms. Thus, to elucidate structure, we will mainly adopt the
perspective of differential forms.

Further, we base our considerations on recent discoveries linking symplectic algebra
and self-adjoint extensions of symmetric operators, see \cite{EVM05} for a survey.
In the context of ordinary differential equations, this connection was intensively
studied by Markus and Everitt during the past few years \cite{EVM99}. They also
extended their investigations to partial differential operators like $\Delta$
\cite{EVM03}. We are going to apply these powerful tools to the special case
of $\curl$ operators. Here, the crucial symplectic space is a Hilbert space
of 1-forms on $\partial\dom$ equipped with the pairing
\begin{gather*}
  [\omega,\eta]_{\partial\dom} := \int\nolimits_{\partial\dom}\omega\wedge\eta\;.
\end{gather*}

We find out, that it is the Hodge decomposition of the trace space for 1-forms on
$\dom$ that allows a classification of self-adjoint extensions of $\curl$: the main
distinction is between boundary conditions that impose closed and co-closed traces
Moreover, further constraints are necessary in the form of vanishing circulation
along certain fundamental cycles of $\partial\dom$. This emerges from an analysis
of the space of harmonic 1-forms on $\partial\dom$ as a finite-dimensional
symplectic space. For all these self-adjoint $\curl$s we show that they possess
a complete orthonormal system of eigenfunctions.

In detail, the outline of the article is as follows: The next section reviews the
connection between vector analysis and differential forms in 3D and 2D.  Then, in the
third section, we introduce basic concepts of symplectic algebra.  Then we summarize
how those can be used to characterize self-adjoint extensions through complete
Lagrangian subspaces of certain factor spaces. The fourth section applies these
abstract results to trace spaces for 1-forms and the corresponding exterior
derivative, that is, the $\curl$ operator. The following section describes important
complete Lagrangian subspaces spawned by the Hodge decomposition of 1-forms on
surfaces. The role of co-homology spaces comes under scrutiny. In the sixth section
we elaborate concrete boundary conditions for self-adjoint $\curl$ operators induced
by the complete Lagrangian subspaces discussed before. The two final sections examine
the spectral properties of the classes of self-adjoint $\curl$s examined before and
explore their relationships with $\curl\curl$ operators. Frequently used notations
are listed in an appendix.

\section{The $\curl$ operator and differential forms}
\label{sec:curl-operator}

In classical vector analysis the operator $\curl$ is introduced as first order
partial differential operator acting on vector fields with three components. Thus,
given a domain $D\subset\bbR^{3}$ we may formally consider $\curl:
\VC^{\infty}_{0}(\dom)\mapsto\VC^{\infty}_{0}(\dom)$ as an unbounded operator on
$\Ltwov$. Integration by parts according to \eqref{eq:curlGreen} shows that this
basic $\curl$ operator is symmetric, hence closable \cite[Ch.~5]{WEI80}. Its closure
is given by the \emph{minimal $curl$ operator}
\begin{gather}
  \label{eq:curlmin}
  \curl_{\min}:\zHcurl\mapsto \Ltwov\;.
\end{gather}
Its adjoint is the \emph{maximal $curl$ operator}, see \cite[Sect.~0]{PIC98a},
\begin{gather}
  \label{eq:curlmax}
  \curl_{\max}:=\curl_{\min}^{\ast} :
  \Hcurl\mapsto \Ltwov\;.
\end{gather}
Note, that $\curl_{\max}$ is no longer symmetric, and neither operator is
self-adjoint. This motivates the search for self-adjoint extensions
$\curl_{s}:\Cd(\curl_{s})\subset\Ltwov\mapsto \Ltwov$ of $\curl_{\min}$. If they
exist, they will satisfy, \textit{c.f.} \cite[Example~1.13]{EVM05},
\begin{gather}
  \label{eq:3}
  \curl_{\min} \subset \curl_{s} \subset \curl_{\max}\;.
\end{gather}
\begin{Remark}
  \label{rem:1}
The classical route in the study of self-adjoint extensions of
symmetric operators is via the famous Stone-von Neumann extension
theory, see \cite[Ch.~6]{WEI80}. It suggests that, after complexification, we examine
the deficiency spaces
\begin{gather}
  \label{eq:11}
  N^{\pm} := \Kern{\curl_{\max}\pm \imath \cdot \mathsf{Id}}\subset
  \Cd(\curl_{\max})\;.
\end{gather}
\begin{lemma}
  \label{lem:DIcurl}
  The deficiency spaces from \eqref{eq:11} satisfy
  $\displaystyle \dim N^{\pm} = \infty$.
\end{lemma}
\begin{proof}
  Let $\VG^{\pm}:\bbR^{3}\setminus\{0\}\mapsto \bbC^{3,3}$ be a fundamental solution
  (dyad) of $\curl\pm \imath$, that is, $\VG = (\curl\mp i)(-1 -
  \nablabf^{T}\nablabf)\VI\Phi$, where $\Phi(\Bx) = \exp(-|\Bx|)/(4\pi|\Bx|)$ is the
  fundamental solution of $-\Delta+1$, and $\nablabf := (\frac{\partial}{\partial
    x_{1}},\frac{\partial}{\partial x_{2}},\frac{\partial}{\partial x_{3}})$.
  Then, for any $\varphibf\in {\VC^{\infty}(\bbR^{3})}_{|\partial\dom}$,
  \begin{gather*}
    \Vu(\Bx) := \int\nolimits_{\partial\dom}\VG(\Bx-\By)\cdot
    \varphibf(\By)\,\mathrm{d}S(\By)\;,\quad \Bx\in \dom\;,
  \end{gather*}
  satisfies $\Vu\in\Hcurl$ and $\curl\Vu\pm \imath \Vu = 0 $.
\end{proof}

From Lemma~\ref{lem:DIcurl} we learn that $N^{\pm}$ reveal little about the structure
governing self-adjoint extensions of $\curl$. Yet the relationship of $\curl$ and
differential forms suggests that there is rich structure underlying self-adjoint
extensions of $\curl_{\min}$.
\end{Remark}

\subsection{Differential forms}
\label{sec:differential-forms}

The $\curl$ operator owes its significance to its close link with the exterior
derivative operator in the calculus of differential forms. We briefly recall its
basic notions and denote by $M$ an $m$-dimensional compact orientable manifold with
boundary $\partial M$. If $M$ is of class $C^{1}$ it can be endowed with a space of
differential forms of degree $k$, $0\leq k\leq m$:
\begin{definition}[Differential $k$-form]
  \label{def:differential-forms}
  A differential form of degree $k$ (in short, a $k$-form) and class $C^{l}$,
  $l\in\mathbb{N}_{0}$, is a $C^{l}$-mapping assigning to each $x\in M$ an
  alternating $k$-multilinear form on the tangent space $T_{x}(M)$. We write
  $\Lambda^{k,l}(M)$ for the vector space of $k$-forms of class $C^{l}$ on $M$,
  $\Lambda^{k,l}(M)=\{0\}$ for $k<0$ or $k>m$.
\end{definition}

Below, we will usually drop the smoothness index $l$, tacitly assuming that
the forms are ``sufficiently smooth'' to allow the respective operations.

The alternating exterior product of multilinear forms gives rise to the exterior
product $\wedge:\Lambda^{k}(M)\times\Lambda^{j}(M)\mapsto\Lambda^{k+j}(M)$ by
pointwise definition. We note the graded commutativity rule $\omega\wedge\eta =
(-1)^{kj}\eta\wedge\omega$ for $\omega\in\Lambda^{k}$, $\eta\in\Lambda^{j}$. Further,
on any piecewise smooth orientable $k$-dimensional sub-manifold of $M$ we can
evaluate the \emph{integral} $\int_{\Sigma}\omega$ of a $k$-form $\omega$ over a
$k$-dimensional sub-manifold $\Sigma$ of $M$ \cite[Sect.~4]{CAR74}.

This connects to the integral view of $k$-forms as entities that describe additive
and continuous (w.r.t. to a suitable deformation topology) mappings from orientable
sub-manifolds of $M$ into the real numbers. This generalized differential forms are
sometimes called currents and are studied in geometric integration theory
\cite{FED69,DRA84}. From this point of view differential forms also make sense
for non-smooth manifolds.

From the integral perspective the transformation (pullback) $\Phi^{*}\omega$ of a
$k$-form under a sufficiently smooth mapping $\Phi:\widehat{M}\mapsto M$ appears
natural: $\Phi^{*}\omega$ is a $k$-form on $\widehat{M}$ that fulfills
\begin{gather}
  \label{eq:6}
  \int\nolimits_{\widehat{\Sigma}}\Phi^{*}\omega =
  \int\nolimits_{\Sigma}\omega
\end{gather}
for all $k$-dimensional orientable sub-manifolds $\widehat{\Sigma}$ of $\widehat{M}$.
We remark that pullbacks commute with the exterior product.

If $i:\partial M\mapsto M$ stands for the inclusion map, then the natural trace of a $k$-form
$\omega\in\Lambda^{k}(M)$ on $\partial M$ is defined as $i^{*}\omega$.  It satisfies
the following commutation relations
\begin{equation}
\label{commutative-relation}
i^*(\omega\wedge\eta)=(i^*\omega)\wedge{}(i^*\eta)
\quad\hbox{ and }\quad\sd{}(i^*\omega)\;=\;i^*(\sd{}\omega)\;.
\end{equation}
The key operation on differential form is the \emph{exterior derivative}
\begin{equation}
  \label{eq:2}
  \sd:\Lambda^{k}(M)\rightarrow\Lambda^{k+1}(M)\;,
\end{equation}
which is connected with integration through the generalized Stokes theorem
\begin{gather}
  \label{eq:4}
  \int\nolimits_{\Sigma}\sd\omega = \int\nolimits_{\partial\Sigma}\omega\quad
  \forall \omega\in \Lambda^{k,0}(M)
\end{gather}
and all orientable piecewise smooth sub-manifolds of $M$. In fact, \eqref{eq:4} can
be used to \emph{define} the exterior derivative in the context of geometric
integration theory. This has the benefit of dispensing with any smoothness
requirement stipulated by the classical definition of $\sd$. We have $\sd^{2}=0$ and,
obviously, \eqref{eq:4} and \eqref{eq:2} imply $\Phi^{*}\circ\sd = \sd\circ\Phi^{*}$.

Since one has the graded Leibnitz formula
\begin{equation}
\sd(\omega\wedge \eta)\;=\;\sd\omega\wedge \eta + (-1)^{k}
  \omega\wedge \sd \eta\quad\forall
  \omega\in\Lambda^{k}(M),\,\eta\in\Lambda^{j}(M)
\end{equation}
exterior derivative and exterior product enter the crucial integration by
parts formula
\begin{gather}
  \label{eq:5}
  \int\nolimits_{\Sigma} \sd\omega\wedge \eta + (-1)^{k}\int\nolimits_{\Sigma}
  \omega\wedge \sd \eta = \int\nolimits_{\partial\Sigma}i^{\ast}\omega\wedge
  i^{\ast}\eta
  \quad\forall
  \omega\in\Lambda^{k}(M),\,\eta\in\Lambda^{j}(M)\;,
\end{gather}
where $\Sigma$ is an orientable sub-manifold of $M$ with dimension $k+j+1$
and canonical inclusion $i:\Sigma\mapsto\partial\Sigma$.

\subsection{Metric concepts}
\label{sec:metric-concepts}

A metric $g$ defined on the manifold $M$ permits us to introduce the \emph{Hodge
  operator} $\hop_{g}:\Lambda^{k}(M)\mapsto \Lambda^{m-k}(M)$. It gives rise to the
inner product on $\Lambda^{k}(M)$
\begin{gather}
  \label{eq:7}
  (\omega,\eta)_{k,M}\;:=\;\int_{M}\omega\wedge*_{g}\eta\;,\quad
  \omega,\eta\in\Lambda^{k}(M)\;.
\end{gather}
Thus, we obtain an $L^{2}$-type norm $\N{\cdot}$ on $\Lambda^{k}(M)$. Completion of
smooth $k$-forms with respect to this norm yields the Hilbert space
$L^{2}(\Lambda^{k}(M))$ of square integrable (w.r.t. $g$) $k$-forms on $M$. Its
elements are equivalence classes of $k$-forms defined almost everywhere on $M$.
Since Lipschitz manifolds possess a tangent space almost everywhere, for them
$L^{2}(\Lambda^{k}(M))$ remains meaningful.  As straightforward is the introduction
of ``Sobolev spaces'' of differential forms, see \cite[Sect.~1]{AFW06},
\begin{gather}
  \label{eq:14}
  W^{k}(\sd,M) := \{ \omega \in L^{2}(\Lambda^{k}(M)):\;
  \sd\omega \in L^{2}(\Lambda^{k}(M))\}\;,
\end{gather}
which are Hilbert spaces with the graph norm. The completion of the subset of
smooth $k$-forms with compact support in $W^{k}(\sd,M)$ is denoted by
$W_{0}^{k}(\sd,M)$.

By construction, the Hodge star operator satisfies
\begin{gather}
  \label{eq:8}
  \hop\hop = (-1)^{(m-k)k}\;\sId\;.
\end{gather}
Now, let us assume $\partial M =\emptyset$. Based on the inner product \eqref{eq:7}
we can introduce the adjoint $\sd^{\ast}:W^{k+1}(\sd,M)\mapsto W^{k}(\sd,M)$ of
the exterior derivative operator by
\begin{equation}
  \label{eq:10}
  (\sd\omega,\eta)_{k+1,M}\;=\;(\omega,\sd^\ast\eta)_{k,M}
  \quad
  \forall \omega\in W^{k+1}(\sd,M),\,\eta\in W^{k}_{0}(\sd,M)\;,
\end{equation}
and an explicit calculation shows that
\begin{equation}
  \label{eq:9}
  \sd^\ast=(-1)^{(mk+1)}\hop\sd\hop:\Lambda^{k+1}\rightarrow\Lambda^{k}\;.
\end{equation}
Furthermore, $\sd^2=0$ implies $(\sd^\ast)^2=0$. Eventually, the Laplace-Beltrami
operator is defined as
\begin{equation} \Delta_{M}\;=\;\sd\sd^\ast+\sd^\ast\sd:\Lambda^k\longrightarrow\Lambda^k.
\end{equation}

\subsection{Vector proxies}
\label{sec:vector-proxies}

Let us zero in on the three-dimensional ``manifold'' $\dom$.  Choosing bases for the
spaces of alternating $k$-multilinear forms, differential $k$-forms can be identified
with vector fields with $\binom{3}{k}$ components, their so-called ``vector proxies''
\cite[Sect.~1]{AFW06}. The usual association of ``Euclidean vector proxies'' in
three-dimensional space is summarized in Table~\ref{tab:vp}. The terminology honours
the fact that the Hodge operators $\hop:\Lambda^{1}(D)\mapsto\Lambda^{2}(D)$ and
$\Lambda^{0}(D)\mapsto\Lambda^{3}(D)$ connected with the Euclidean metric of 3-space
leave the vector proxies invariant (this is not true in 2D since $\hop^2=-1$ on
1-forms).  In addition the exterior product of forms is converted into the pointwise
Euclidean inner product of vector fields. Thus, the inner product
$(\cdot,\cdot)_{k,\dom}$ of $k$-forms on $D$ becomes the conventional $L^{2}(\dom)$
inner product of the vector proxies.  Further, the spaces $W^{k}(\sd,\dom)$ boil down
to the standard Sobolev spaces $H^{1}(\dom)$ (for $k=0$), $\Hcurl$ (for $k=1$),
$\Hdiv$ (for $k=2$), and $\Ltwo$ (for $k=3$).

\begin{table}[!htb]
  \centering
\begin{tabular}{ll}
  \hline\hline
  Differential form $\omega$ & Related function $u$ or vector field $\Vu$ \\\hline
  $\Bx\mapsto{}\omega(\Bx)$ & $u(\Bx):=\omega(\Bx)$ \\[3pt]
  $\Bx\mapsto{}\{\textbf{v}\mapsto\omega(\Bx)(\textbf{v})\}$ & $\textbf{u}(\Bx)\cdot\textbf{v}:=\omega(\Bx)(\textbf{v})$ \\[3pt]
  $\Bx\mapsto{}\{(\textbf{v}_1,\textbf{v}_2)\mapsto\omega(\Bx)(\textbf{v}_1,\textbf{v}_2)\}$ & $\textbf{u}(\Bx)\cdot\big(\textbf{v}_1\times{\textbf{v}_2}\big):=\omega(\Bx)(\textbf{v}_1,\textbf{v}_2)$ \\[3pt]
  $\Bx\mapsto{}\{(\textbf{v}_1,\textbf{v}_2,\textbf{v}_3)\mapsto\omega(\Bx)(\textbf{v}_1,\textbf{v}_2,\textbf{v}_3)\}$ & $u(\textbf{x})\;\hbox{det}(\textbf{v}_1,\textbf{v}_2,\textbf{v}_3):=\omega(\Bx)(\textbf{v}_1,\textbf{v}_2,\textbf{v}_3)$ \\[3pt]\hline
  \hline
\end{tabular}
\caption{The standard choice of vector proxy $u$,$\Vu$ for a differential form
  $\omega$ in $\bbR^3$. Here, $\cdot$ denotes the Euclidean inner product of
  vectors in $\bbR^{3}$, whereas $\times$ designates the cross product.}
\label{tab:vp}
\end{table}

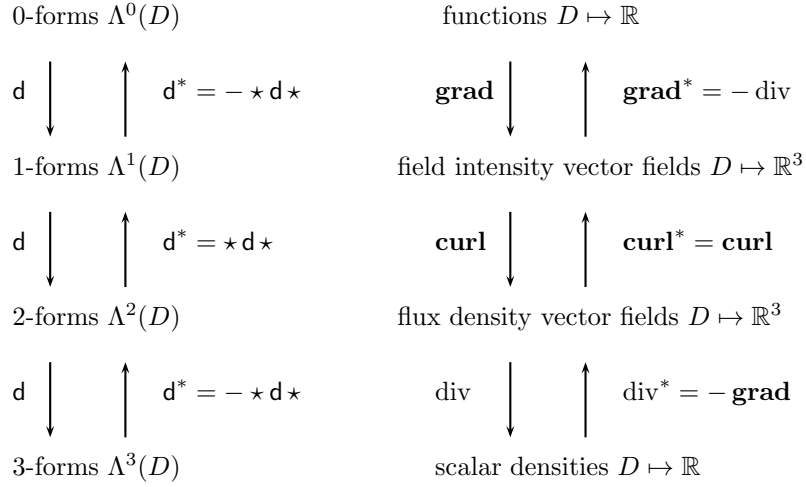
\begin{figure}[!htb]
  \centering
  \setlength{\unitlength}{1cm}
  {
\begin{picture}(4,6.5)
  \psline{<-}(1.5,.5)(1.5,1.5)
  \psline{<-}(1.5,2.5)(1.5,3.5)
  \psline{<-}(1.5,4.5)(1.5,5.5)
  \psline{->}(2.5,.5)(2.5,1.5)
  \psline{->}(2.5,2.5)(2.5,3.5)
  \psline{->}(2.5,4.5)(2.5,5.5)
  \put(1,0){3-forms $\Lambda^{3}(D)$}
  \put(1,2){2-forms $\Lambda^{2}(D)$}
  \put(1,4){1-forms $\Lambda^{1}(D)$}
  \put(1,6){0-forms $\Lambda^{0}(D)$}
  \put(1,5){$\sd$}
  \put(1,3){$\sd$}
  \put(1,1){$\sd$}
  \put(3,5){$\sd^*=-\hop\sd\hop$}
  \put(3,3){$\sd^*=\hop\sd\hop$}
  \put(3,1){$\sd^*=-\hop\sd\hop$}
\end{picture}\hspace{2cm}
\begin{picture}(4,6.5)
  \psline{<-}(1.5,.5)(1.5,1.5)
  \psline{<-}(1.5,2.5)(1.5,3.5)
  \psline{<-}(1.5,4.5)(1.5,5.5)
  \psline{->}(2.5,.5)(2.5,1.5)
  \psline{->}(2.5,2.5)(2.5,3.5)
  \psline{->}(2.5,4.5)(2.5,5.5)
  \put(0.5,0){scalar densities $D\mapsto\bbR$}
  \put(0,2){flux density vector fields $D\mapsto\bbR^{3}$}
  \put(0,4){field intensity vector fields $D\mapsto\bbR^{3}$}
  \put(0.6,6){functions $D\mapsto\bbR$}
  \put(.5,5){$\grad$}
  \put(.5,3){$\curl$}
  \put(.5,1){$\Div$}
  \put(3,5){$\grad^\ast=-\Div$}
  \put(3,3){$\curl^\ast=\curl$}
  \put(3,1){$\Div^\ast=-\grad$}
\end{picture} }
  \caption{Differential operators and their relationship with exterior derivatives}
  \label{fig:diffop}
\end{figure}

Using Euclidean vector proxies, the $\curl$ operator turns out to be an incarnation
of the \emph{exterior derivative} for 1-forms. More generally, the key first order
differential operators of vector analysis arise from exterior derivative operators,
see Figure~\ref{fig:diffop}. Please note that, since the Hodge operator is
invisible on the vector proxy side, $\curl$ can as well stand for the operator
\begin{gather}
  \label{eq:1}
  \curl\quad \longleftrightarrow\quad
  \hop \sd : \Lambda^{1}(\dom)\mapsto \Lambda^{1}(\dom)\;,
\end{gather}
which is naturally viewed as an unbounded operator on $L^{2}(\Lambda^{1}(\dom))$.
Thus, \eqref{eq:1} puts the formal $\curl$ operator introduced above in the framework
of differential forms on $\dom$.

Translated into the language of differential forms, the Green's formula
\eqref{eq:curlGreen} becomes a special version of \eqref{eq:5} for $k=j=1$.  However,
due to \eqref{eq:8}, \eqref{eq:curlGreen} can also be stated as
\begin{gather}
  \label{eq:IPcurl}
  (\hop\sd\omega,\eta)_{1,\dom} - (\omega,\hop\sd\eta)_{1,\dom} =
  \int\nolimits_{\partial\dom} i^{\ast}\omega \wedge i^{\ast}\eta\;,\quad
  \omega,\eta \in W^{1}(\sd,\dom)\;.
\end{gather}

A metric on $\bbR^{3}$ induces a metric on the embedded 2-dimensional manifold
$\partial \dom$. Thus, the Euclidean inner product on local tangent spaces
becomes a meaningful concept and Euclidean vector proxies for $k$-forms
on $\partial \dom$, $k=0,1,2$, can be defined as in Table~\ref{tab:vp},
see Table~\ref{tab:2dvp}.

\begin{table}[!htb]
  \centering
  \begin{tabular}{ll}
  \hline\hline
  Differential forms & Related function $u$ or vector field $\Vu$\\\hline
  $\Bx\mapsto{}\omega(\Bx)$ & $u(\Bx):=\omega(\Bx)$\\[3pt]
  $\Bx\mapsto{}\{\textbf{v}\mapsto\omega(\Bx)(\textbf{v})\}$ & $\textbf{u}(\Bx)\cdot\textbf{v}:=\omega(\Bx)(\textbf{v})$ \\[3pt]
  $\Bx\mapsto{}\{(\textbf{v}_1,\textbf{v}_2)
  \mapsto\omega(\Bx)(\textbf{v}_1,\textbf{v}_2)\}$ &
  $u(\Bx)\;\hbox{det}(\textbf{v}_1,\textbf{v}_2,\Vn(\Bx)):=
  \omega(\Bx)(\textbf{v}_1,\textbf{v}_2)$ \\[3pt]\hline
  \hline
\end{tabular}
  \caption{Euclidean vector proxies for differential forms on $\partial\dom$. Note
    that the test vectors $\Vv,\Vv_{1},\Vv_{2}$ have to be chosen from the tangent space
    $T_{\Bx}(\partial \dom)$.}
  \label{tab:2dvp}
\end{table}

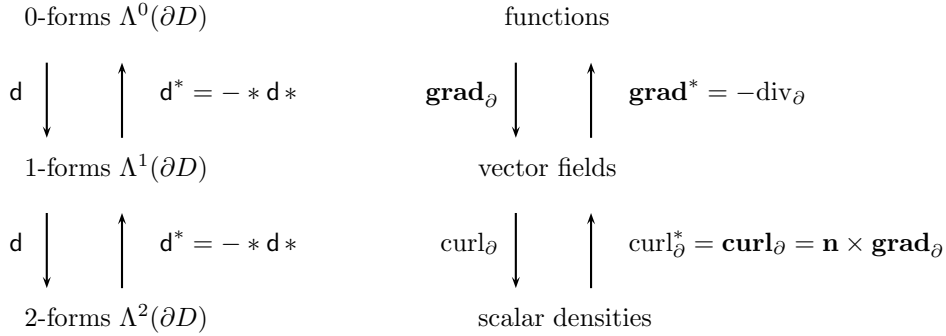
\begin{figure}[!htb]
  \centering
  \setlength{\unitlength}{1cm} {
    \begin{picture}(4,4.5)(2,0)
      \psline{<-}(1.5,.5)(1.5,1.5) \psline{<-}(1.5,2.5)(1.5,3.5)
      \psline{->}(2.5,.5)(2.5,1.5) \psline{->}(2.5,2.5)(2.5,3.5)
      \put(1.2,0){2-forms $\Lambda^{2}(\partial\dom)$}
      \put(1.2,2){1-forms $\Lambda^{1}(\partial\dom)$}
      \put(1.2,4){0-forms $\Lambda^{0}(\partial\dom)$}
      \put(1,3){$\sd$}
      \put(1,1){$\sd$}
      \put(3,3){$\sd^*=-*\sd*$}
      \put(3,1){$\sd^*=-*\sd*$}
    \end{picture}
    \hspace{2cm}
    \begin{picture}(4,4.5)(2,0)
      \psline{<-}(1.5,.5)(1.5,1.5)
      \psline{<-}(1.5,2.5)(1.5,3.5)
      \psline{->}(2.5,.5)(2.5,1.5)
      \psline{->}(2.5,2.5)(2.5,3.5)
      \put(1,0){scalar densities}
      \put(1,2){vector fields}
      \put(1.35,4){functions}
      \put(.3,3){$\bgrad$}
      \put(.5,1){$\bscurl$} \put(3,3){$\grad^{\ast}=-\bDiv$}
      \put(3,1){$\bscurl^{\ast}=\bcurl=\mathbf{n}\times\bgrad$}
    \end{picture} }
  \caption{Exterior derivative for Euclidean vector proxies on 2-manifolds}
  \label{fig:vp2d}
\end{figure}

This choice of vector proxies leads to convenient vector analytic expressions
for the trace operator $i^{\ast}$:
\begin{gather*}
  \left\{
    \begin{array}{lllclll}
      \omega\in\Lambda^0(D):&{}i^*\omega
      &\longleftrightarrow&\gamma u(\Bx)&:=&u(\Bx),\;&u:\dom\mapsto\bbR\;,
      \\
      \omega\in\Lambda^1(D):&{}i^*\omega
      &\longleftrightarrow&\gamma_t \Vu(\Bx)&:=&\Vu(\Bx)-(\Vu(\Bx)\cdot\Vn(\Bx))
      \Vn(\Bx),&\Vu:D\mapsto\bbR^{3}\;,
      \\
      \omega\in\Lambda^2(D):&{}i^*\omega
      &\longleftrightarrow&\gamma_n \Vu(\Bx)&:=&\Vu(\Bx)\cdot\Vn(\Bx),&
      \Vu:D\mapsto\bbR^{3}\;,
      \\
      \omega\in\Lambda^3(D):&{}i^*\omega
      &\longleftrightarrow&0\;,
    \end{array}
  \right.
\end{gather*}
where $\Bx\in\partial\dom$. Further, the customary vector analytic surface
differential operators realize the exterior derivative for vector proxies,
see Figure~\ref{fig:vp2d}.

\begin{Remark}
  \label{rem:vp}
  Vector proxies offer an isomorphic model for the calculus of differential forms.
  However, one must be aware that the choice of bases and, therefore, the description
  of a differential form by a vector proxy, is somewhat arbitrary.  In particular, a
  change of metric of space suggests a different choice of vector proxies for which the
  Hodge operators reduce to the identity. Thus, metric and topological aspects are
  hard to disentangle from a vector analysis point of view. This made us prefer
  the differential forms point of view in the remainder of the article.
\end{Remark}

\section{Self-adjoint extensions and Lagrangian subspaces}
\label{sec:self-adjo-extens}

First, we would like to recall some definitions of symplectic geometry. Then, we will
build a symplectic space associated to a closed symmetric operator. The reader can
refer to \cite{EVM03,EVM05} for a more detailed treatment.

\subsection{Concepts from symplectic geometry}
\label{sec:basic-defin-prop}

Symplectic geometry offers an abstract framework to deal with self-adjoint
extensions of symmetric operators in Hilbert spaces. Here we briefly review
some results. More information is available from \cite{MDS95}.

\begin{definition}[Symplectic space]
  \label{def:31}
  A real symplectic space ${S}$ is a real linear space equipped with a
  symplectic pairing $[\cdot,\cdot]$ (symplectic bilinear form, symplectic product)
  \begin{equation}
    \left\{
      \begin{array}{ll}
        [\cdot,\cdot]:S\times S\longrightarrow\bbR, \\[6pt]
        [\alpha_1u_1+\alpha_2u_2,v]\;=\;\alpha_1[u_1,v]\;+\;\alpha_2[u_2,v],
        & (linearity)\\[6pt]
        [u,v]\;=\;-{[v,u]}, & \text{(skew symmetry)}\\[6pt]
        [u,S]\;=\;0\quad\Longrightarrow\quad{}u=0 & \text{(non-degeneracy)}
      \end{array}
    \right.
  \end{equation}
\end{definition}
\begin{definition}
  \label{def:32}
  Let $L$ be a linear subspace of the symplectic space $S$
\begin{itemize}
\item[(i)] The \emph{symplectic orthogonal}  of $L$ is $L^\symporth\;=\;\{u\in S:\;[u,L]=0\}$;
\item[(ii)] $L$ is a \emph{Lagrangian} subspace, if $L \subset L^\symporth$ i.\ e. $[u,v]=0$
  for all $u$ and $v$ in $L$;
\item[(iii)] A \emph{Lagrangian} subspace $L$ is complete, if $L^{\symporth}=L$.
\end{itemize}
\end{definition}

In the case of finite dimensional symplectic spaces, symplectic bases offer
a convenient way to build complete Lagrangian subspaces, see \cite[Example~2]{EVM99}.
\begin{definition}
  Let $(S,[\cdot,\cdot])$ be a real symplectic space with dimension $2n$ (the
  dimension has to be even so that the pairing $[\cdot,\cdot]$ can be non
  degenerate). A \emph{symplectic basis} $\{u_{i}\}_{i=1}^{2n}$ of $S$ is a basis of $S$
  satisfying
\begin{equation}
    [u_i,u_j]=\VJ_{i,j}\quad\hbox{ with }\mathbf{J}=
  \begin{bmatrix}
0_{n\times{}n}&\bf{I}_{n\times{}n}\\-\bf{I}_{n\times{}n}&0_{n\times{}n}
  \end{bmatrix}
\end{equation}
\end{definition}
 Simple linear algebra proves the existence of such bases:
\begin{lemma}
  For any symplectic space with finite dimension $2n$, there exists a (non unique) symplectic basis.
\end{lemma}

\begin{Remark}
  \label{rem:ls}
  As soon as we have found a symplectic basis $\{u_i\}_{i=1}^{2n}$,
  it provides many complete Lagrangian subspaces
  \begin{itemize}
  \item the $n$ first vectors $\{u_i\}_{i=1}^{n}$ of a symplectic
    basis span a complete Lagrangian  subspace.
  \item the $n$ last vectors $\{u_i\}_{i=n+1}^{2n}$ of a symplectic
    basis span a complete Lagrangian  subspace.
  \item for any $\sigma:[\![1,n]\!]\mapsto[\![0,1]\!]$,
    $\{u_{i+\sigma(i)n}\}_{i=1}^{n}$ is a complete Lagrangian  subspace.
  \end{itemize}
\end{Remark}

We recall some more facts about finite dimensional symplectic spaces 
\begin{lemma}
  \label{lem:fds}
  Every complete Lagrangian subspace of a finite dimensional
  symplectic space $S$ of dimension $2n$ is $n$-dimensional. Moreover,
  it  possesses a basis that can be extended to a
  symplectic basis of $S$.
\end{lemma}

\subsection{Application to self-adjoint extensions of a symmetric operator}
\label{sec:appl-self-adjo}

Let $H$ be a real Hilbert space and $\sT$ a closed symmetric linear operator with
dense domain $\Cd(\sT)\subset{}H$. We denote by $\sT^*$ its adjoint. Let us first recall,
see \cite{WEI80}, that each self-adjoint extension of $\sT$ is a restriction of
$\sT^*$, which is classically written as
\begin{equation}
  \sT\subset\sT^s\subset\sT^*.
\end{equation}
Next, introduce a degenerate symplectic pairing on $\Cd(\sT^*)$ by
\begin{equation}\label{eq:[.,.]}
  [\cdot,\cdot]:\Cd(\sT^*)\times\Cd(\sT^*)\longrightarrow\bbR
  \quad\hbox{ such that }\quad[\su,\sv]\;=\;(\sT^*\su,\sv)\;-\;(\su,\sT^*\sv).
\end{equation}
From the definition of $\sT^{\ast}$, the symmetry of $\sT$, and the fact
$\sT^{\ast\ast} = \sT$, we infer that, see \cite[Appendix]{EVM99},
\begin{equation}\label{eq:sym_pai}
  \left\{
    \begin{array}{lllllll}
      [u+u_0,v+v_0]\;=\;[\su,\sv],\quad\forall{}u_0,v_0\in\Cd(\sT),\;
      \forall{}u,v\in\Cd(\sT^*),\\[6pt]
      u\in{}\Cd(\sT^*),\;[u,v]\;=\;0,\quad\forall{}v\in{}\Cd(\sT^*),
      \quad\implies\quad{}u\in{}\Cd(\sT).
    \end{array}
  \right.
\end{equation}
As a consequence, we obtain a symplectic factor space, see Appendix
of \cite{EVM99},
\begin{lemma}\label{lem:sym}
  The space $S=\Big(\Cd(\sT^*)/\Cd(\sT),[\cdot,\cdot]\Big)$ is a symplectic space.
\end{lemma}
The graph norm on $\Cd(\sT^*)$ induces a factor norm on $S$ and, due to \eqref{eq:sym_pai},  the symplectic pairing $[\cdot,\cdot]$ is continuous with respect to this norm
\begin{gather*}
  |[\su,\sv]|^{2} \leq \bigl(\N{u}^{2}+\N{\mathsf{T}^{\ast}u}^{2}\bigr)\cdot
  \bigl(\N{v}^{2}+\N{\mathsf{T}^{\ast}v}^{2}\bigr)\quad
  \forall u\in\Cd(\sT^\ast),\,v\in\Cd(\sT^\ast)\;,
\end{gather*}
Let $L\oplus\Cd(\sT)$ denotes the pre-image of $L$ under the factor map
$\Cd(\sT^*)\mapsto S$.
\begin{corollary}
  The symplectic orthogonal $V^{\symporth}$ of any subspace $V$ of $S$ is closed
  (in the factor space topology).
\end{corollary}
Any linear subspace $L$ of $S$ defines an extension $\sT_{L}$ of $\mathsf{T}$ through
\begin{gather}
  \label{eq:16}
  \sT \subset \sT_{L} := {\sT^{\ast}}_{|L\oplus\Cd(\sT)} \subset \sT^{\ast}\;.
\end{gather}
This relationship allows to characterize self-adjoint extensions of $\sT$ by means
of the symplectic properties of the associated subspaces $L$. This statement is made
precise in the Glazman-Krein-Naimark Theorem, see Theorem~1 of
\cite[Appendix]{EVM99}.
\begin{theorem}[Glazman-Krein-Naimark Theorem symplectic version]
  \label{GKN-symplectic}
  The mapping $L\mapsto\sT_{L}$ is a bijection between the space of complete Lagrangian
  subspaces of $S$ and the space of self-adjoint extensions of $\sT$. The inverse
  mapping is given by
  \begin{equation}
    L=\Cd(\sT_{L})/\Cd(\sT)\;.
  \end{equation}
\end{theorem}

\section{Symplectic space for $\curl$}
\label{sec:sympl-space-curl}

\newcommand*{\SYP}[1]{\left[{#1}\right]_{\partial\dom}} Evidently, the unbounded
$\curl$ operators introduced in Section~\ref{sec:curl-operator} (resorting to the
vector proxy point of view) fits the framework of the preceding section and
Theorem~\ref{GKN-symplectic} can be applied. To begin with, from \eqref{eq:curlmin}
and \eqref{eq:curlmax} we arrive at the symplectic space
\begin{gather}
  \label{eq:12}
  S_{\curl} := \Hcurl/\zbHcurl\;.
\end{gather}
By \eqref{eq:curlGreen} it can be equipped with a symplectic pairing that can
formally be written as
\begin{gather}
  \label{eq:18}
  \SYP{\su,\sv} := \int\nolimits_{\partial \dom}
  (\su(\By)\times\sv(\By))\cdot\Vn(\By)\,\mathrm{d}S(\By)\;,
\end{gather}
for any representatives of the equivalence classes of $S_{\curl}$. From
\eqref{eq:12} it is immediate that $S$ is algebraically and topologically
isomorphic to the natural trace space of $\Hcurl$.

By now this trace space is well understood, see the seminal work of Paquet
\cite{PAC82} and \cite{BUC99,BUC99a,BUF02,BCS00} for the extension to generic
Lipschitz domains. To begin with, the topology of $S_{\curl}$ is intrinsic,
that is, with $\dom':=\bbR^{3}\setminus \overline{\dom}$, the norm of 
\begin{gather}
  \label{eq:17}
  S_{\curl}^{c} := \Hcurl[\dom']/\zbHcurl[\dom']
\end{gather}
is equivalent to that of $S_{\curl}$; both spaces are isomorphic algebraically and
topologically. This can be proved appealing to an extension theorem for $\Hcurl$.
Moreover, the pairing $\SYP{\cdot,\cdot}$ identifies $S_{\curl}$ with its dual
$S_{\curl}'$:
\begin{lemma}
  The mapping $S_{\curl}\mapsto S_{\curl}'$, $\su\mapsto\{\sv\mapsto \SYP{\su,\sv}\}$ is
  an isomorphism.
\end{lemma}
\begin{proof}
  Given $\su\in S_{\curl}$, let $\Vu\in\Hcurl$ solve
  \begin{gather}
    \label{eq:19}
    \curl\curl\Vu+\Vu = 0 \quad\text{in }\dom\;,\quad
    \gamma_{t}\Vu = \su\quad\text{on }\partial\dom\;.
  \end{gather}
  Set $\Vv:=\curl\Vu\in\Hcurl$ and $\sv := \gamma_{t}\Vv\in S_{\curl}$. By
  \eqref{eq:curlGreen}
  \begin{eqnarray*}
    \SYP{\su,\sv} &=& \int\nolimits_{\dom} 
    \curl\Vu\cdot\Vv - \curl\Vv\cdot\Vu\,\mathrm{d}\Vx \\
    &=& \int\nolimits_{\dom} |\curl\Vu|^{2} + |\Vu|^{2}\,\mathrm{d}\Vx =
    \NHcurl{\Vv}\NHcurl{\Vu} \geq \N{\sv}_{S_{\curl}}\N{\su}_{S_{\curl}}\;,
  \end{eqnarray*}
  as $\NHcurl{\Vv}=\NHcurl{\Vu}$. We immediately conclude
  \begin{gather*}
    \sup\limits_{\sv\in S_{\curl}}\frac{|[\su,\sv]|}{\N{\sv}_{S_{\curl}}} \geq \N{\su}_{S_{\curl}}\;.
  \end{gather*}
\end{proof}

\renewcommand{\Hhv}{\textbf{H}_{t}^{\frac{1}{2}}(\partial\dom)}
\newcommand{\Hhmv}{\textbf{H}_{t}^{-\frac{1}{2}}(\partial\dom)}
\newcommand{\Hmhc}{\textbf{H}^{-\frac{1}{2}}(\bscurl,\partial\dom)}
The trace space also allows a characterization via surface differential
operators. It relies on the space $\Hhv$ of tangential surface traces
of vector fields in $(\Hone)^{3}$ and its dual $\Hhmv := (\Hhv)'$.
Then one finds that, algebraically and topologically, $S_{\curl}$
is isomorphic to
\begin{gather}
  \label{eq:22}
  S_{\curl} \cong \Hmhc := \{\sv\in \Hhmv:\; \bscurl\sv\in \Hm[\partial\dom]{-\frac{1}{2}}\}\;.
\end{gather}
The intricate details and the proper definition of $\bscurl$ can be 
found in \cite{BCS00}.

When we adopt the perspective of differential forms, the domain of
$\curl_{\max}$ is the Sobolev space $H^{1}(\sd,\dom)$ of 1-forms. Thus,
$S_{\curl}$ has to be viewed as a trace space of 1-forms, that is, 
a space of 1-forms (more precisely, 1-currents) on $\partial\dom$.
In analogy to \eqref{eq:14} and \eqref{eq:22} we adopt the notation
\begin{gather}
  \label{eq:23}
  S_{\curl} \cong \tscurl\;.
\end{gather}
Please observe, that the corresponding symbol for the trace space of
$H^{0}(\sd,\dom)$ will be $\tsgrad$ (and not $W^{\frac{1}{2},0}(\sd,\dom)$ 
as readers accustomed to the conventions used with Sobolev spaces might expect).

In light of \eqref{eq:IPcurl}, the symplectic pairing on $\tscurl$ can be expressed as
\begin{gather}
  \label{eq:21}
  \SYP{\omega,\eta} := \int\nolimits_{\partial\dom} \omega\wedge\eta\;,\quad
  \omega,\eta\in\tscurl\;.
\end{gather}
Whenever, $\tscurl$ is treated as a real symplectic space, the pairing \eqref{eq:21}
is assumed. The most important observation about \eqref{eq:21} is that the pairing
$[\cdot,\cdot]$ is utterly metric-free!

Now we can specialize Theorem~\ref{GKN-symplectic} to the $\curl$ operator. We give
two equivalent versions, one for Euclidean vector proxies, the second  for 1-forms:

\begin{theorem}[GKN-theorem for $\curl$, vector proxy version]
  \label{thm:GKNproxies}
  The mapping which associates to $L\subset \Hmhc$ the $\curl$ operator with domain 
  \begin{gather*}
    \Cd(\curl_L) := \{\Vv\in\Hcurl:\; \gamma_{t}(\Vv) \in L\}
  \end{gather*}
  is a bijection between the set of complete Lagrangian subspaces of $\Hmhc$ and the self-adjoint extensions of $\curl_{\min}$.
\end{theorem}

\begin{theorem}[GKN-theorem for $\curl$, version for 1-forms]
  \label{thm:GKN}
  The mapping which associates to $L\subset \tscurl$ the ${\hop\sd}$ operator with domain
  \begin{gather*}
    \Cd({\hop\sd}_L) := \{\eta\in W^{1}(\sd,\dom):\; i^{\ast}\eta \in L\}
  \end{gather*}
  is a bijection between the set of complete Lagrangian subspaces of $\tscurl$ and
  the self-adjoint extensions of $\hop\sd$ defined on $W^{1}_{0}(\sd,\dom)$.
\end{theorem}

We point out that the constraint $i^{\ast}\eta \in L$ on traces amounts to imposing
linear \emph{boundary conditions}. In other words, the above theorems tell us, that
self-adjoint extensions of $\curl_{\min}$ will be characterized by demanding
particular boundary conditions for their argument vector fields, \textit{cf.}
\cite{PIC98a}.

\begin{Remark}
  \label{rem:mi}
  Thanks to \eqref{commutative-relation} the symplectic pairing on $\tscurl$ commutes
  with the pullback. Thus, if $\dom,\widehat{\dom}\subset\bbR^{3}$ are connected by a
  Lipschitz homomorphism $\Phi:\widehat{\dom}\mapsto\dom$, we find that
  $\Phi^{\ast}:\Lambda^{1}(\partial\dom)\mapsto\Lambda^{1}(\partial\widehat{\dom})$
  provides a bijective mapping between the complete Lagrangian subspaces of $\tscurl$
  and of $\tscurl[\partial\widehat{\dom}]$. Thus, pulling back the domain of a
  self-adjoint extension of $\curl_{\min}$ on $\dom$ to $\widehat{\dom}$ will give a
  valid domain for a self-adjoint extension of $\curl_{\min}$ on $\widehat{\dom}$.  In
  short, self-adjoint extensions of $\curl_{\min}$ are invariant under bijective
  continuous transformations. This is a very special feature of $\curl$, not shared,
  for instance, by the Laplacian $-\Delta$.
\end{Remark}
 
\section{Hodge theory and consequences}
\label{sec:hodge-theory-cons}

Now we study particular subspaces of the trace space $\tscurl$. We will take for
granted a metric on $\partial\dom$ that induces a Hodge operator $\hop$.

\subsection{The Hodge decomposition}
\label{sec:hodge-decomposition}

Let us first recall the well-known Hodge decomposition of spaces of square-integrable
differential 1-forms on $\partial\dom$. For a more general exposition we refer to
\cite{MOR66}:

\begin{lemma}
  \label{lem:hd}
  We have the following decomposition, which is orthogonal w.r.t. the
  inner product of $L^{2}(\Lambda^{1}(\partial\dom))$:
  \begin{gather*}
    L^{2}(\Lambda^{1}(\partial\dom) = d W^{0}(\sd,\partial\dom) \oplus \hop\sd
    W^{0}(\sd,\partial\dom) \oplus \Ch^{1}(\partial\dom)\;.
  \end{gather*}
\end{lemma}

Here, $\Ch^{1}(\partial\dom)$ designates the finite-dimensional space of harmonic 1-forms on
$\partial\dom$:
\begin{gather}
  \label{eq:13}
  \Ch^{1}(\partial\dom) := \{\omega\in L^{2}(\Lambda^{1}(\partial\dom)):\;
  \sd\omega=0\;\text{and}\; \sd\hop\omega=0\}\;.
\end{gather}

In terms of Euclidean vector proxies, the space $L^{2}(\Lambda^{1}(\partial\dom)$ is
modelled by the space $\mathbf{L}^{2}_{t}(\partial\dom)$ of square integrable
tangential vector fields on $\partial\dom$. Then, the decomposition of
Lemma~\ref{lem:hd} reads
\begin{gather*}
  \mathbf{L}^{2}_{t}(\partial\dom) = \bgrad H^{1}(\partial\dom) \oplus
  \bcurl H^{1}(\partial\dom) \oplus \Ch^{1}(\partial\dom)\;,\\
  \Ch^{1}(\partial\dom) := \{\Vv\in\mathbf{L}^{2}_{t}(\partial\dom):\;
  \bcurl\Vv=0\;\text{and}\; \bDiv\Vv=0\}\;.
\end{gather*}
The Hodge decomposition can be extended to $\tscurl$ on Lipschitz domains, as was
demonstrated in \cite[Sect.~5]{BCS00} and \cite{BUF00}. There the authors showed
that, with a suitable extension of the surface differential operators, that
\begin{gather}
  \label{eq:15}
  \Hmhc = \bgrad \Hh[\partial\dom] \oplus \bcurl \Hm[\partial\dom]{\frac{3}{2}}
  \oplus \Ch^{1}(\partial\dom)\;,
\end{gather}
where, formally,
\begin{gather}
  \label{eq:20}
  \Hm[\partial\dom]{\frac{3}{2}} := \Delta_{\partial\dom}^{-1}
  H_{\ast}^{-\frac{1}{2}}(\partial\dom)\;,\quad H_{\ast}^{-\frac{1}{2}}(\partial\dom)
  := \{v\in H^{-\frac{1}{2}}(\partial\dom):\,
  \int\nolimits_{\partial\dom_i}v\,\mathrm{d}S = 0 \}\;.
\end{gather}
with $\partial\dom_i$ the connected components of $\partial{}\dom$.\\[6pt]
For $C^{1}$-boundaries this space agrees with the trace space of $H^{2}(\dom)$.

The result \eqref{eq:15} can be rephrased in the calculus of differential forms:
\begin{theorem}[Hodge decomposition of trace space]
  \label{Hodge-forms}
  We have the following orthogonal decomposition
  \begin{gather}
    \tscurl = \sd \tsgrad \oplus \hop\sd W^{\frac{3}{2},0}(\partial\dom)
    \oplus  \Ch^{1}(\partial\dom)\;,
  \end{gather}
  whith
  \begin{gather}
    W^{\frac{3}{2},0}(\partial\dom) :=
    \Delta_{\partial\dom}^{-1}\bigl\{
    \varphi\in{}W^{-\frac{1}{2},0}(\partial{}D\;:\;\langle\varphi,\mathbf{1}\rangle_{\partial{}D_i} = 0\bigr\} \;.
  \end{gather}
with $\partial\dom_i$ the connected components of $\partial{}\dom$.
\end{theorem}

The first subspace in the decomposition of Theorem~\ref{Hodge-forms} comprises
only closed 1-forms, because
\begin{gather}
  \label{eq:24}
  \sd\bigl(\sd \tsgrad \bigr) = 0 \;.
\end{gather}
The second subspace contains only co-closed 1-forms, since
\begin{gather}
  \label{eq:25}
  \sd^{\ast}\bigl(\hop\sd W^{\frac{3}{2},0}  \bigr) = 0 \;.
\end{gather}

The Hodge decomposition hinges on the choice of the Hodge operator
$\hop$. Consequently, it depends on the underlying metric on $\partial\dom$.

\subsection{Lagrangian properties of the Hodge decomposition}
\label{sec:lagr-prop-hodge}

We find that the subspaces contributing to the Hodge decomposition of
Theorem~\ref{Hodge-forms} can be used a building blocks for (complete) Lagrangian
subspaces of $\tscurl$.
\begin{proposition}
  \label{pro:clo}
  The linear space $\sd \tsgrad$ is a Lagrangian subspace of $\tscurl$
  (w.r.t. the symplectic pairing $\SYP{\cdot,\cdot}$)
\end{proposition}
\begin{proof}
  We have to show that
  \begin{equation}
    \SYP{\sd\omega,\sd\eta}=0
    \quad\forall\omega,\eta\in{}\tsgrad\;.
  \end{equation}
  By density, we need merely consider $\omega$, $\eta$ in $W^{0}(\sd,\partial\dom)$.
  In this case, it is immediate from Stokes' Theorem ($\partial\dom$ has no
  boundary)
  \begin{equation}
    [\sd\omega,\sd\eta]_{\partial{}D}\;=\;
    \int_{\partial{}D}\sd\omega\wedge\sd\eta\;=\;\int_{\partial{}D}\omega\wedge\sd^2\eta=0\;.
  \end{equation}
\end{proof}
\begin{proposition}
  \label{pro:coclo}
  The linear space $\hop\sd{}W^{\frac{3}{2},0}(\partial{}D)$ is a
  Lagrangian subspace of $\tscurl$.
\end{proposition}
\begin{proof}
  The proof is the same as above, except that one has to use that $\hop$ is an
  isometry with respect to the inner product induced by it:
  \begin{equation}
    [\hop\sd\omega,\hop\sd\eta]_{\partial{}D}
    =
    \int_{\partial{}D}\hop\sd\omega\wedge\hop\sd\eta=
    - \int_{\partial{}D}\sd\omega\wedge\sd\eta=-\int_{\partial{}D}\omega\wedge\sd^2\eta=0\;.
  \end{equation}
\end{proof}

In a similar way we prove the next proposition.
\begin{proposition}
  \label{pro:harm} The space of harmonic 1-forms
  $\Ch^1(\partial{}D)$ is symplectically orthogonal to
  $\sd{}\tsgrad$ and $\hop\sd{}W^{\frac{3}{2},0}(\partial{}D)$.
\end{proposition}

The Hodge decomposition of Theorem~\ref{Hodge-forms} offers a tool for the evaluation
of the symplectic pairing ${[\cdot,\cdot]}_{\partial\dom}$.  Below, tag the three
components of the Hodge decomposition of Theorem~\ref{Hodge-forms} by subscripts $0$,
$\perp$ and $H$: for $\omega,\eta\in\tscurl$
\begin{equation}
  \label{hodge-decom}
  \omega\;=\;\sd\omega_0\;+\;\hop\sd\omega_{\perp}\;+\;\omega_H\quad\hbox{ and }
  \quad\eta\;=\;\sd\eta_0\;+\;\hop\sd\eta_{\perp}\;+\;\eta_H\;.
\end{equation}
Note that the forms $\omega_0$ and $\omega_{\perp}$ are not unique since the kernels of
$\sd$ and $\hop\sd$ are not empty (they contain the piecewise constants on conected
components of $\partial\dom$).

Taking into account the symplectic orthogonalities stated in Propositions \ref{pro:clo},
\ref{pro:coclo}, and \ref{pro:harm}
\begin{equation}
  \left\{
    \begin{array}{cccccccccc}
      [\sd\omega_0,\sd\eta_0]_{\partial{}D}&=&[\hop\sd\omega_{\perp},\hop\sd\eta_{\perp}]_{\partial{}D}&=&[\sd\omega_0,\eta_H]_{\partial{}D}&=&[\hop\sd\omega_{\perp},\eta_H]_{\partial{}D}\\
      &=&[\omega_H,\sd\eta_0]_{\partial{}D}&=&[\omega_H,\hop\sd\eta_{\perp}]_{\partial{}D}&=&0
    \end{array}
  \right.
\end{equation}
we see that we can compute the symplectic pairing on $\tscurl$ according to
\begin{equation}\label{comp:symp}
  [\omega,\eta]_{\partial{}D}\;=\;[\sd\omega_0,\hop\sd\eta_{\perp}]_{\partial{}D}+
  [\hop\sd\omega_{\perp},\sd\eta_0]_{\partial{}D}\;+\;[\omega_H,\eta_H]_{\partial{}D}\;.
\end{equation}
It can also been expressed in terms of the $L^{2}$-inner product (more precisely, its extension to
duality pairing) as
\begin{equation}
  [\omega,\eta]_{\partial{}D}\;=\;(\sd\omega_0,\sd\eta_{\perp})_{1,\partial{}D} -
  (\sd\omega_{\perp},\sd\eta_0)_{1,\partial{}D}\;+\;[\omega_H,\eta_H]_{\partial{}D}\;.
\end{equation}

\subsection{The symplectic space $\Ch^{1}(\partial\dom)$}
\label{sec:sympl-space-ch1p}

Let us recall that the space of harmonic 1-forms on $\partial\dom$ (a 2 dimensional
compact $C^\infty$-manifold without boundary) is a finite dimensional linear space
with
\begin{equation}
  \dim(\Ch^1(\partial\dom))\;=\;2g,
\end{equation}
with $g$ the genus of the boundary, that is, the first Betti number of $\dom$. The
reader can refer to Theorem 5.1, Proposition 5.3.1 of \cite{BOT82} and Theorem 7.4.3
of \cite{MOR66}.

Since the set of harmonic vector fields is stable with respect to the Hodge operator
(note  that $\hop\hop = -1$ for 1-forms on $\partial\dom$)
\begin{equation}
  \left\{
    \begin{array}{cccccccc}
      \eta\in\Ch^1(\partial\dom)&\quad\Longrightarrow\quad&
      \eta\in L^{2}(\Lambda^{1}(\partial\dom)),&\sd\eta=0,&\sd\hop\eta=0\\[4pt]
      &\quad\Longrightarrow\quad&
      \hop\eta\in L^{2}(\Lambda^{1}(\partial\dom)),&\sd\hop(\hop\eta)=0,&\sd(\hop\eta)=0\\[4pt]
      &\quad\Longrightarrow\quad&\hop\eta\in\Ch^1(\partial\dom),
    \end{array}
  \right.
\end{equation}
Thus we find that the pairing $\SYP{\cdot,\cdot}$ is non-degenerate on $\Ch^1(\partial\dom)$:
\begin{equation}
  \Big([\omega_H,\eta_H]_{\partial\dom}\;=\;0,\quad\forall\eta_H\in\Ch^1(D)\Big)\Longrightarrow[\omega_H,\hop\omega_H]_{\partial\dom} = (\omega_{H},\omega_{H})_{1,\partial\dom} = 0 \;.
\end{equation}
\begin{lemma}\label{lem:finitesymp}
  The space of harmonic 1-forms $\Ch^1(\partial\dom)$ is a symplectic space with
  finite dimension when equipped with the symplectic pairing
  $[\cdot,\cdot]_{\partial\dom}$. It is a finite-dimensional symplectic subspace of $\tscurl$.
\end{lemma}

\section{Some examples of self-adjoint $\curl$ operators}
\label{sec:some-examples-self}

Starting from the Hodge decomposition of Theorem~\ref{Hodge-forms}, we now identify
important classes of self-adjoint extensions of $\curl$. We rely on a generic Riemannian metric
on $\partial\dom$ and the associated Hodge operator.

\subsection{Self-adjoint $\curl$ associated with closed traces}
\label{sec:self-adjoint-curl}

In this section we aim to characterize the complete Lagrangian subspaces $L$ of
$W^{-\frac{1}{2},1}(\sd,\partial{}D)$ (equipped with $[\cdot,\cdot]_{\partial{}D}$) which
contain only closed forms:
\begin{equation}
  \label{eq:26}
  L\subset Z^{-\frac{1}{2},1}(\partial{}D):=\{\eta\in\tscurl:\;\sd\eta=0\}\;.
\end{equation}
Hodge theory (see Theorem \ref{Hodge-forms}) provides the tools to
study these Lagrangian subspaces, since we have the following result:
\begin{lemma}
  \label{Hodge-closed}
  The set of closed 1-forms in $\tscurl$ admits the following direct (orthogonal)
  decomposition
  \begin{equation}\label{Hodgeclosed}
    Z^{-\frac{1}{2},1}(\partial{}D)\;=\;\sd{}W^{-\frac{1}{2},0}(\sd,\partial{}D)\;
    \oplus\;\Ch^1(\partial{}D)\;.
  \end{equation}
\end{lemma}
\begin{proof}
  For $\omega\in
  Z^{-\frac{1}{2},1}(\partial{}D)$, $\hop\sd\omega_{\perp}$ of \eqref{hodge-decom} satisfies
  \begin{equation}
    \sd(\hop\sd\omega_{\perp})=0,\;\sd\hop(\hop\sd\omega_{\perp})=0,\;(\hop\sd\omega_{\perp})_H=0
  \end{equation}
  which implies that $\hop\sd\omega_{\perp}=0$ and yields the assertion of the lemma.
\end{proof}

The next result is important, as it states a necessary condition for
the existence of Lagrangian subspaces included in
$Z^{-\frac{1}{2},1}(\partial{}D)$.
\begin{lemma}
  \label{lem:clo:lag}
  The space $Z^{-\frac{1}{2},1}(\partial{}D)$ includes its
  symplectic orthogonal $\sd{}W^{-\frac{1}{2},0}(\sd,\partial{}D)$.
\end{lemma}
\begin{proof}
  Recall from Definition~\ref{def:32} that the symplectic orthogonal of
  $Z^{-\frac{1}{2},1}(\partial{}D)$ is defined as the set
  \begin{equation}
    \{\omega\in{}W^{-\frac{1}{2},1}(\sd,\partial{}D):\;
    [\omega,\eta]_{\partial{}D}=0,\quad\forall{}\eta\in{}\sZ^{-\frac{1}{2},1}(\partial{}D)\}\;.
  \end{equation}
  Using Theorem \ref{Hodge-forms} for $\omega=\sd\omega_0+\hop\sd\omega_{\perp}+\omega_H$
  and Lemma \ref{Hodge-closed} for $\eta=\sd\eta_0+\eta_H$, we have with \eqref{comp:symp}:
  \begin{equation}
    \label{hodge-decom-symp}
    [\omega,\eta]_{\partial{}D}\;=\;
    [\hop\sd\omega_{\perp},\sd\eta_0]_{\partial{}D}+[\omega_H,\eta_H]_{\partial{}D}.
  \end{equation}
  This implies with $\eta=\hop\omega_H=\eta_H$
  (here we use the stability of $\Ch^1(\partial{}D)$ with respect to the Hodge operator)
  \begin{equation}
    [\omega,\hop\omega_H]_{\partial{}D}=[\omega_H,\hop\omega_H]_{\partial{}D}
    =\int_{\partial{}D}\omega_H\wedge\hop\omega_H=0\quad\implies\omega_H=0\;,
  \end{equation}
  and, for $\eta=\sd\eta_0$ with $\eta_0=\omega_{\perp}\in{}W^{3/2,0}(\partial{}D)$
  \begin{equation}
    [\omega,\sd\omega_{\perp}]_{\partial{}D}\;=\;[\hop\sd\omega_{\perp},\sd\omega_{\perp}]_{\partial{}D}
    \;=\;-\int_{\partial{}D}\sd\omega_{\perp}\wedge\hop\sd\omega_{\perp}\quad\Longrightarrow\quad\sd\omega_{\perp}=0
    \;.
  \end{equation}
  Hence, we have $\omega=\sd\omega_0$ (and $\omega_H=0$).  The converse holds due to
  \eqref{hodge-decom-symp}.
\end{proof}

Lemma~\ref{lem:clo:lag} tells us that, when restricted to the subspace of closed
forms, the bilinear pairing $[\cdot,\cdot]_{\partial{}D}$ becomes degenerate.  More
precisely, on the subset of closed forms, one can use the splitting
\eqref{hodge-decom-symp} and evaluate $[\cdot,\cdot]_{\partial{}D}$ on
$Z^{-\frac{1}{2},1}(\partial\dom)$ according to
\begin{equation}
  \label{coho-harm-symp}
  [\omega,\eta]_{\partial{}D}\;=\;
  [\omega_H,\eta_H]_{\partial{}D}\;,\quad
  \forall{}\omega,\;\eta\in{}\sZ^{-\frac{1}{2},1}(\partial{}D)\;.
\end{equation}
Hence, this pairing depends only on the harmonic components. Thus another message of
Lemma~\ref{lem:clo:lag} is that $[\cdot,\cdot]_{\partial{}D}$ furnishes a
well-defined non-degenerate symplectic pairing, when considered on the co-homology
factor space
\begin{equation}
  \mathbb{H}^1(\partial{}D,\bbR)\;=\;
  Z^{-\frac{1}{2},1}(\partial{}D))/\sd{}W^{-\frac{1}{2},0}(\sd,\partial{}D)\;.
\end{equation}
This space is algebraically, topologically and symplectically isomorphic to
$\Ch^1(\partial{}D)$, the space of harmonic 1-forms, see
Section~\ref{sec:sympl-space-ch1p}.

This means that all the complete Lagrangian subspaces $L$ of $\tscurl$ contained in
$Z^{-\frac{1}{2},1}(\partial{}D)$ are related to complete Lagrangian subspaces
$L_{\Ch}$ of $\Ch^1(\partial{}D)$ (or equivalently to the complete Lagrangian
subspace $L_{\mathbb{H}}$ of $\mathbb{H}^1(\partial{}D,\bbR)$) by
\begin{equation}
  \label{rel:lagrclo}
  L\;=\;\sd{}W^{-\frac{1}{2},0}(\sd,\partial{}D)\oplus L_\Ch\quad\quad\hbox{ (or, equivalently, }
  L_\Ch\;=\;L/\sd{}W^{-\frac{1}{2},0}(\sd,\partial{}D)\hbox{)}\;.
\end{equation}
Thus, we have proved the following lemma (the symplectic pairing
$[\cdot,\cdot]_{\partial{}D}$ is used throughout)
\begin{lemma}
  \label{theo:lag_closed}
  There is a one-to-one correspondance between the complete
  Lagrangian subspaces $L$ of the symplectic space $\tscurl$ satisfying
  \begin{equation}
    L\subset{}Z^{-\frac{1}{2},1}(\partial{}D)
  \end{equation}
  and the complete Lagrangian subspaces $L_\Ch$ of $\Ch^1(\partial{}D)$. The
  bijection is given by \eqref{rel:lagrclo}.
\end{lemma}

Theorem \ref{thm:GKN} and Lemma \ref{theo:lag_closed} lead
to the characterization of the self-adjoint $\curl$ operators whose
domains contain only functions with closed traces.
\begin{theorem}
  \label{theo:closed}
  There is a one-to-one corresondance between the set of all selfadjoint
  $\curl$ operators $\hop\sd_S$ satisfying
  \begin{equation}
    \label{eq:29}
    \Cd(\hop\sd_S)\;\subset\;
    \Big\{\omega\in{}W^{1}(\sd,\dom)\;\Big|\;i^{\ast}\omega
    \in{}Z^{-\frac{1}{2},1}(\partial{}D)\Big\}
  \end{equation}
  and the set of complete Lagrangian subspaces $L_{\Ch}$ of
  $\Ch^{1}(\partial{}\dom)$.
  They are related according to
  \begin{equation}
    \label{eq:28}
    \Cd(\hop\sd_S)\;=\;\Big\{\omega\in{}W^{1}(\sd,\dom)\;
    \Big|\;i^{\ast}\omega\in{}\sd{}W^{-\frac{1}{2},0}(\sd,\partial{}D)\oplus L_\Ch\Big\}\;.
  \end{equation}
\end{theorem}
Obviously, the constraint
\begin{equation}
  i^{\ast}\omega\in{}\sd{}W^{-\frac{1}{2},0}(\sd,\partial{}D)\oplus L_\Ch
\end{equation}
is a boundary condition, since it involves only the boundary of the
domain $D$. In addition, we point out that no metric concepts enter
in \eqref{eq:28}, \textit{cf.} Section~\ref{sec:sympl-space-curl}.

\begin{Remark}
  \label{rem:closed}
  Now, assume the domain $D$ to feature trivial topology, that is, the genus of $D$ is zero,
  and the space of harmonic forms is trivial. Theorem \ref{theo:closed} reveals that
  there is only one self-adjoint $\hop\sd$ with domain containing only forms
  with closed traces
  \begin{equation}
    \Cd(\hop\sd_S)\;=\;\Big\{\omega\in{}W^1(\sd,\dom)\;\Big|\;\sd{}(i^*\omega)=0\Big\}.
  \end{equation}
  In terms of vector proxies, this leads to the self-adjoint $\curl$ operator with domain
  \begin{equation}
    \Cd(\curl_S)\;=\;\Big\{{\su}\in{}W^0(\curl,D)\;\Big|\;
    \curl({\su})\cdot\Vn=0\hbox{ on }\partial{}D\Big\}\;,
  \end{equation}
  which has been investigated in \cite{PIC98,YOG90}.  In case $\dom$ has non-trivial
  topology, then $\hbox{dim}(\Ch^1(\partial{}D))=2g\neq 0$ and one has to examine the
  complete Lagrangian subspaces of $\Ch^1(\partial{}D)$, which is postponed to
  Section~\ref{sec:compl-lagr-ch1p}.
\end{Remark}

\subsection{Self-adjoint $\curl$ based on co-closed traces}
\label{sec:self-curl-assoc}

In this section we seek to characterize those Lagrangian subspaces $L$ of $\tscurl$ that
contain only co-closed forms, \textit{ie.}
\begin{equation}\label{def:coclosed}
  L\subset{}\Big\{\omega\in{}W^{-1/2,1}(\sd,\partial{}D):\;
  \mathsf{d}\hop\omega=0\Big\}\;.
\end{equation}
The developments are parallel to those of the previous
section, because, as is illustrated by \eqref{comp:symp}, from a symplectic point
of view, the subspaces of closed and co-closed 1-forms occuring in the Hodge
decomposition of Theorem~\ref{Hodge-forms}, are symmetric. For the sake of
completeness, we give the details, nevertheless.
\begin{lemma}
  \label{Hodge-coclosed}
  The subspace of co-closed 1-forms of $W^{-1/2,1}(\sd,\partial{}D)$ admits the
  following orthogonal decomposition
\begin{equation}
  \label{Hodgecoclosed}
  \Big\{\omega\in{}W^{-1/2,1}(\sd,\partial{}D):\;\mathsf{d}\hop\omega=0\Big\}
  \;=\;\hop\sd{}W^{3/2,0}(\partial{}D)\;\oplus\;\Ch^1(\partial{}D)\;.
\end{equation}
\end{lemma}
\begin{proof}
  For $\omega$ co-closed, $\sd\omega_0$ in \eqref{hodge-decom} satisfies
  \begin{equation}
    \sd(\sd\omega_0)=0,\;\sd\hop(\sd\omega_0)=0,\;(\sd\omega_0)_H=0
  \end{equation}
  which implies that $\sd\omega_0=0$ and proves \eqref{Hodgecoclosed}.
\end{proof}

The next result is important as it states a necessary condition for
the existence of Lagrangian subspaces comprising only co-closed
forms.
\begin{lemma}
  \label{lem:coclo:lag}
  The symplectic orthogonal of the subspace of
  co-closed forms of $W^{-1/2,1}(\sd,\partial{}D)$
  is $\hop\sd{}W^{3/2,0}(\partial{}D)$.
\end{lemma}
\begin{proof}
  Use the definition of the symplectic orthogonal as
  \begin{equation}
    \bigl\{\omega\in{}W^{-1/2,1}(\sd,\partial{}D):\;
    [\omega,\eta]_{\partial{}D}=0\quad
    \forall{}\eta\hbox{ co-closed}\bigr\}\;.
  \end{equation}
  Using Theorem \ref{Hodge-forms} for $\omega=\sd\omega_0+
  \hop\sd\omega_{\perp}+\omega_H$ and Lemma \ref{Hodge-coclosed} for
  $\eta=\hop\sd\eta_{\perp}+\eta_H$, \eqref{comp:symp} gives
  \begin{equation}
    \label{hodge-decom-symp-coc}
    [\omega,\eta]_{\partial{}D}\;=\;
    [\sd\omega_0,\hop\sd\eta_{\perp}]_{\partial{}D}+[\omega_H,\eta_H]_{\partial{}D}\;.
  \end{equation}
  Choosing $\eta=\hop\omega_H=\eta_H$
  (here we use the stability of $\Ch^1(\partial{}D)$ with respect to the Hodge
  operator) this implies
  \begin{gather*}
    [\omega,\hop\omega_H]_{\partial{}D}
    =[\omega_H,\hop\omega_H]_{\partial{}D}=
    \int_{\partial{}D}\omega_H\wedge\hop\omega_H=0\quad\implies\omega_H=0,
  \end{gather*}
  and, for $\eta=\hop\sd\eta_{\perp}$ with
  $\eta_{\perp}\in{}W^{\frac{3}{2},0}(\partial{}D)$
  \begin{gather*}
    [\omega,\hop\sd\eta_{\perp}]_{\partial{}D}\;
    =\;[\sd\omega_0,\hop\sd\eta_{\perp}]_{\partial{}D}\;=\;
    \int_{\partial{}D}\sd\omega_0\wedge\hop\sd\eta_{\perp}\;=\;0
    \quad\Longrightarrow\quad\sd\hop(\sd\omega_0)=0\;.
  \end{gather*}
  Moreover, one has $\sd(\sd{}\omega_0)=0$ and $(\sd\omega_0)_H=0$,
  which shows that $\sd\omega_0=0$.

  Hence, we have $\omega=\hop\sd\omega_{\perp}$ ($\sd{}\omega_0=0$ and $\omega_H=0$).
  The other inclusion holds due to \eqref{hodge-decom-symp-coc}.
\end{proof}

\begin{Remark}
  Formally, in the proof we have used $\eta=\hop\sd\eta_\perp$ with $\eta_\perp=\omega_0$ which
  shows that
  \begin{equation}
    [\omega,\eta]_{\partial{}D}\;=\;
    [\sd\omega_0,\hop\sd\omega_0]_{\partial{}D}\;=
    \int_{\partial{}D}\sd\omega_0\wedge\hop\sd\omega_0\quad\Longrightarrow\quad\sd\omega_0=0.
  \end{equation}
  However, the lack of regularity of $\omega_0$ does not allow this straightforward
  computation.
\end{Remark}

When restricted to the space of co-closed forms, the bilinear pairing
$[\cdot,\cdot]_{\partial{}D}$ becomes degenerate.  However, due to Lemma
\ref{lem:coclo:lag}, it is a non-degenerate symplectic product on the
co-homology factor space.
\begin{equation}
  \Big\{\omega\in{}W^{-1/2,1}(\sd,\partial{}D):\;\sd\hop
  \omega=0\Big\}\Big/\hop\sd{}W^{3/2,0}(\partial{}D)\;,
\end{equation}
which can be identified with $\Ch^1(\partial{}D)$.  Indeed, also on the subset of
co-closed forms, one can evaluate $[\cdot,\cdot]_{\partial{}D}$ by means of
\eqref{coho-harm-symp}.
\begin{lemma}
  \label{theo:lag_coclosed}
  The complete Lagrangian subspaces $L$ of $W^{-1/2,1}(\sd,\partial{}D)$ containing
  only co-closed forms are one-to-one related to the complete Lagrangian subspaces
  $L_\Ch$ of $\Ch^1(\partial{}D)$ by
  \begin{equation}\label{rel:lagr}
    L\;=\;\hop\sd{}W^{3/2,0}(\partial{}D)\oplus L_\Ch\;.
\end{equation}
\end{lemma}
Theorem \ref{thm:GKN} and Lemma \ref{theo:lag_coclosed} lead
to the characterization of the self-adjoint $\curl$ operators based
on coclosed forms:
\begin{theorem}
  \label{theo:coclosed}
  There is a one to one correspondance between the set of all selfadjoint operators
  $\hop\sd_S$ satisfying
  \begin{equation}
    \Cd(\hop\sd_S)\;\subset\;
    \Big\{\omega\in{}W^1(\sd,\partial{}D):\;\sd\hop(i^*\omega)=0\Big\}
  \end{equation}
  and the set of complete Lagrangian subspaces $L_\Ch$ of $\Ch^{1}(\partial{}D)$
  equipped with $[\cdot,\cdot]_{\partial{}D}$. The underlying bijection is
  \begin{equation}
    \label{eq:30}
    \Cd(\hop\sd_S)\;=
    \;\Big\{\omega\in{}W^{1}(\sd,\dom):\;i^*\omega
    \in{}\hop\sd{}W^{3/2,0}(\partial{}D)\oplus L_\Ch\Big\}\;.
  \end{equation}
\end{theorem}
\begin{Remark}
Let $D$ be a domain with trivial topology. Then there is only one self-adjoint
operator $\hop\sd$ whith domain containing only forms whose traces are coclosed
\begin{equation}
  \Cd(\hop\sd_S)\;=\;\Big\{\omega\in{}W^0(\sd,\Omega):\;\sd{}\hop(i^*\omega)=0\Big\}.
\end{equation}
In terms of Euclidean vector proxies, we obtain the self-adjoint $\curl$ operator
with domain
\begin{equation}
  \Cd(\curl_S)\;=\;\Big\{{\Vu}\in{}\Hcurl:\;
  \bDiv({\gamma_t(\su)})=0\hbox{ on }\partial{}D\Big\}.
\end{equation}
On the contrary, if $D$ has non trivial topology, then one has to identify the
complete Lagrangian subspaces of $\Ch^1(\partial{}D)$. This is the topic of the
next section.
\end{Remark}

\subsection{Complete Lagrangian subspaces of  $\Ch^{1}(\partial\dom)$}
\label{sec:compl-lagr-ch1p}

The goal is to give a rather concrete description of the boundary conditions implied
by \eqref{eq:28} and \eqref{eq:30}. Concepts from topolgy will be pivotal.

To begin with we exploit a consequence of the long Mayer-Vietoris exact sequence in
co-homology \cite{BOT82}, namely the algebraic isomorphisms \cite{KOT90a}
\begin{gather}
  \label{eq:31}
  \Ch^{1}(\partial\dom) \cong
  \mathbb{H}^{1}(\partial\dom;\bbR) \cong
  i^{\ast}_{\mathrm{in}}\mathbb{H}^{1}(\dom;\bbR) +
  i^{\ast}_{\mathrm{out}}\mathbb{H}^{1}(\dom';\bbR)\;.
\end{gather}
Here, $\mathbb{H}^{1}(\dom)$ is the co-homology space $Z^{1}(\sd,\dom)/\sd
W^{0}(\sd,\dom)$, and $i_{\mathrm{in}}:\partial\dom\mapsto\dom$,
$i_{\mathrm{out}}:\partial\dom\mapsto \dom'$ stand for the canonical inclusion maps.
We also point out that {\cite{KOT90a}}
\begin{gather}
  \label{eq:35}
  \tfrac{1}{2}\dim \mathbb{H}^{1}(\partial\dom;\bbR) = \dim\mathbb{H}^{1}(\dom;\bbR) = \dim
  \mathbb{H}^{1}(\dom';\bbR) = g\;,
\end{gather}
where $g\in\bbN_{0}$ is the genus of $\dom$.

Next, we find bases of $\mathbb{H}^{1}(\dom;\bbR)$ and
$\mathbb{H}^{1}(\dom';\bbR)$ using the Poincar\'e duality between
co-homology spaces and relative homology spaces%
\footnote{In this
  article, we denote by 
  \begin{itemize} \item $\bbH_i(A;R)$ the $i^{th}$
    homology group of $A$ with coefficients in $R$;
  \item $\bbH^i(A;R)$ the $i^{th}$
    co-homology space of $A$ with coefficients in $R$;
  \item $\bbH_i(A,B;R)$ the $i^{th}$ relative homology group of $A$ relative to $B$ with
    coefficients in $R$;
  \item $\bbH^i(A,B;R)$ the $i^{th}$ relative
    co-homology space of $A$ relative to $B$ with coefficients in $R$.
  \end{itemize}}

\begin{gather}
  \label{eq:34}
  \mathbb{H}^{1}(\dom;\bbR) \cong \mathbb{H}_{2}(\dom,\partial\dom;\bbR)\;.
\end{gather}
 Consider the relative homology groups (with coefficients in $\mathbb{Z}$)
\begin{equation}
  \mathbb{H}_2(D,\partial{}D;\bbZ)\quad\hbox{ and }\quad{}
  \mathbb{H}_2(D',\partial{}D;\bbZ)
\end{equation}
as integer lattices in the vector spaces
\begin{equation}
  \mathbb{H}_2(D,\partial{}D;\bbR)\quad\hbox{ and }\quad{}
  \mathbb{H}_2(D',\partial{}D;\bbR)\;.
\end{equation}
In \cite{KOT87}, it is shown that these lattices as Abelian groups
are torsion free, and that representatives of homology classes can
be realized as orientable embedded surfaces. More precisely, we can
find $2g$ compact orientable embedded Seifert surfaces, or "cuts"
\begin{equation}
  S_i,\;S'_i\;,\quad{}1\leqslant{}i\leqslant{}g\;,
\end{equation}
such that their equivalence classes under appropriate homology
relations form bases for the following associated
lattices\footnote{Throughout the paper $\langle \cdot\rangle$ denotes the
operation of taking the (relative) homology class of a cycle}
\begin{equation}
  \{\langle S_i\rangle\}_{i=1}^{g}\quad{}\hbox{ for }\mathbb{H}_2(D,\partial{}D;\bbZ),\quad
  \{\langle S'_i\rangle\}_{i=1}^{g}\quad{}\hbox{ for }\mathbb{H}_2(D',\partial{}D;\bbZ)
\end{equation}
and
\begin{equation}
  \{\langle \partial{}S_i'\rangle,\langle \partial{}S_i\rangle\}_{i=1}^{g}\quad\hbox{ for }
  \mathbb{H}_1(\partial{}D;\bbZ).
\end{equation}
In other words, the boundaries $\partial{}S_i,\partial{}S_i'$ provide fundamental
non-bounding cycles on $\partial\dom$.

In \cite{KOT90a}, it was established that the set of surfaces
$\{\langle S_i\rangle\}_{i=1}^{g}\cup\{\langle S_i'\rangle\}_{i=1}^{g}$ can be chosen so that
they are ``dual to each other''. Here, this duality is expressed
through the intersection numbers of their boundaries, see Chapter 5
of \cite{GRK04}.

\newcommand{\Int}{\operatorname{Int}}
\begin{lemma}
  \label{lem:linking}
  The set of surfaces $\{\langle S_i\rangle\}_{i=1}^{g}\cup\{\langle S_i'\rangle\}_{i=1}^{g}$ can be chosen
  such that the intersection pairing on $\mathbb{H}_1(\partial{}D;\bbZ)$ can be
  reduced to ($1\leq i,j\leq g$)
  \begin{equation}
    \label{inter_pairing}
    \left\{
      \begin{array}{l}
        \Int(\langle \partial{}S_i\rangle,\langle \partial{}S'_j\rangle)=\delta_{i,j}\;,\\
        \Int(\langle \partial{}S'_i\rangle,\langle \partial{}S_j\rangle)=-\delta_{i,j}\;.
      \end{array}
    \right.
\end{equation}
\end{lemma}

Furthermore, when the boundaries of these surfaces are "pushed out''
of their respective regions of definition, we get curves in the
complementary region
\begin{equation}\label{bound}
  \partial{}S'_i\longrightarrow{}C_i,\quad\partial{}S_i\longrightarrow{}C_i'.
\end{equation}
The homology classes of these curves form bases for homology
lattices as follows
\begin{equation}
  \{\langle C_i\rangle\}_{i=1}^{g}\quad\hbox{ for }\mathbb{H}_1(D;\bbZ),\quad\hbox{ and }\quad
  \{\langle C'_i\rangle\}_{i=1}^{g}\quad\hbox{ for }\mathbb{H}_1(D';\bbZ).
\end{equation}
This paves the way for a construction of bases of the co-homology
spaces on $\dom$ and $\dom'$ \cite{KOT90a}:
\begin{lemma}
  \label{lem:h1bas}
  The co-homology classes generated by the closed 1-form in the sets
  defined for $1\leq{}i\leq{}g$
  \begin{align*}
    & \Big\{\kappa_{i}\in L^{2}(\Lambda^{1}(\dom)):\;\sd\kappa_{i}=0
    \text{ and }\; \int\nolimits_{C_{j}}\kappa_{i} = \delta_{ij}\hbox{ for }1\leq j\leq g\Big\} \\
    & \Big\{\kappa_{i}'\in L^{2}(\Lambda^{1}(\dom')):\;\sd\kappa_{i}=0
    \text{ and } \int\nolimits_{C_{j}'}\kappa_{i}' = \delta_{ij}\hbox{ for }1\leq j\leq g\Big\}
  \end{align*}
  form bases of $\mathbb{H}^1(D;\bbR)$  and $\mathbb{H}^{1}(D';\bbR)$, respectively.
\end{lemma}

For instance, $\kappa_{i}$ can be obtained as the piecewise exterior
derivative of a 0-form on $\dom\setminus S_{i}$ that has a jump of
height 1 across $S_{i}$. An analoguous statement holds for
$\kappa_{i}'$ with $S_{i}$ replaced with $S_{i}'$. More precisely,
one has for $1\leq{}i\leq{}g$
\begin{equation}\label{kappapsi}
\exists\psi_i\in{}W^{0}(\sd,\dom):\kappa_i=\sd\psi_i\quad\hbox{ on
}D\setminus{}S\quad\hbox{ and }\quad[\psi_i]_{S_j}=\delta_{i,j}
\end{equation}
\begin{equation}\label{kappapsiprime}
\exists\psi'_i\in{}W^{0}(\sd,\dom):\kappa'_i=\sd\psi'_i\quad\hbox{
on }D'\setminus{}S'\quad\hbox{ and
}\quad[\psi'_i]_{S_j'}=\delta_{i,j}
\end{equation}
with $[\cdot]_\Gamma$ denoting the jump across $\Gamma$.

\newcommand{\dps}{\displaystyle}
\begin{lemma}
  \label{lem:3}
  For $1\leqslant{}m,n\leqslant{}g$, we have
  \begin{equation}
    \label{eq:36}
    \begin{array}{llllllllll}
      a)&\dps\int_{\partial{}D}i^*_{in}(\kappa_m)\wedge{}i^*_{in}(\kappa_n)\;=\;0,\\[16pt]
      b)&   \dps \int_{\partial{}D}i^*_{out}(\kappa'_m)\wedge{}i^*_{out}(\kappa'_n)\;=\;0.
    \end{array}
  \end{equation}
\end{lemma}
\begin{proof}
  To establish a) we rewrite the integral as one over $D$, as the following calculation shows
  \begin{equation*}
    \left\{
      \begin{array}{llllllllll}
        \dps\int_{\partial{}D}i^*_{in}(\kappa_m)\wedge{}i^*_{in}(\kappa_n)&=&
        \dps\int_{\partial{}D}i^*_{in}(\kappa_m\wedge{}\kappa_n)\;=\;\int_{D}\sd(\kappa_m\wedge{}\kappa_n)\\[6pt]&=&\dps\int_{D}(\sd\kappa_m)\wedge{}\kappa_n-\kappa_m\wedge{}(\sd\kappa_n)\;=\;0.
      \end{array}
    \right.
  \end{equation*}
  Similarly, b) follows from an analogous calculation where
  $\partial{}D=-\partial{}D'$ with forms defined on $D'$.
\end{proof}

\begin{lemma}
  \label{lem:4}
  For $1\leqslant{}i,\;j\leqslant{}g$, we have
  \begin{equation}
    \int_{\partial{}D}i^*_{in}(\kappa_i)\wedge{}i^*_{out}(\kappa'_j)\;=\;
    \delta_{i,j}.
  \end{equation}
\end{lemma}

\begin{proof}
Let us represent the 1-form $\kappa_i$ by means the 0-form $\psi_i$, which jumps
across $S_{i}$, see \eqref{kappapsi}. Taking into account that
$\sd i^*_{out}\kappa_j'=0$, we get
\begin{equation}
i^*_{in}\kappa_i\wedge{}i^*_{out}\kappa_j'=\sd{}i^*_{in}\psi_i\wedge{}i^*_{out}\kappa_j'
=\sd{}\big(i^*_{in}\psi_i\wedge{}i^*_{out}\kappa'_j\big).
\end{equation}
Applying Stokes Theorem leads to (one has to take care of the
orientation)
\begin{equation}
\dps\int_{\partial{}D}i^*_{in}\kappa_i\wedge{}i^*_{out}\kappa_j'
\;=\;
\dps\int_{\partial{}S_i}\big[i^*_{in}\psi\wedge{}i^*_{out}\kappa'\big]_{\partial{}S_i}
+
\int_{\partial{}S'_j}\big[i^*_{in}\psi\wedge{}i^*_{out}\kappa'\big]_{\partial{}S_i}.
\end{equation}
By \eqref{kappapsi}, we get
\begin{equation}
\dps\int_{\partial{}D}i^*_{in}\kappa_i\wedge{}i^*_{out}\kappa_j'
\;=\;
 \dps\int_{\partial{}S_i}\kappa'_j \;+\;0.
\end{equation}
Since $\partial{}S_i\in\langle C'_i\rangle$, the result follows from \eqref{bound}
and Lemma \ref{lem:h1bas}.
\end{proof}

\begin{Remark}
  When the cuts do not satisfy \eqref{inter_pairing}, a generalization of Lemma
  \ref{lem:4} takes the form
  \begin{equation}
    \int_{\partial{}D}i^*_{in}(\kappa_i)\wedge{}i^*_{out}(\kappa'_j)\;=\;
    \Int(\langle \partial{}S_i\rangle,\langle \partial{}S'_j\rangle)\;.
  \end{equation}
with $\Int(\langle \partial{}S_i\rangle,\langle \partial{}S'_j\rangle)$ the intersection
number of $\langle \partial{}S_i\rangle$ and $\langle \partial{}S'_j\rangle$, see
\cite{GRK04}.
\end{Remark}

Now, take \eqref{inter_pairing} for granted. Write $\kappa_{H,i}$, $\kappa_{H,i}'$,
$1\leq i\leq g$, for the unique harmonic 1-forms, i.e.,
$\kappa_{H,i},\kappa_{H,i}'\in\Ch^{1}(\partial\dom)$, such that
\begin{gather}
  \label{eq:37}
  i_{\mathrm{in}}^{\ast}\kappa_{i} = \kappa_{H,i} + \sd\alpha\quad,\quad
  i_{\mathrm{out}}^{\ast}\kappa_{i}' = \kappa_{H,i}' + \sd\beta\;,
\end{gather}
for some $\alpha,\beta\in L^{2}(\Lambda^{0}(\partial\dom))$. Combining Lemmas
\ref{lem:3} and \ref{lem:4} gives the desired symplectic basis of the
space of harmonic 1-forms on $\partial\dom$:

\begin{lemma}
  \label{lem:harmbas}
  The set $\{\kappa_{H,i},\kappa_{H,i}'\}_{i=1}^{g}$ is a symplectic basis
  of $\Ch^{1}(\partial\dom)$.
\end{lemma}

Obviously, since the trace preserves integrals and integrating a closed form over a
cycle evaluates to zero, the 1-forms $\kappa_{H,i}$ and $\kappa_{H,i}'$ inherit the
integral values over fundamental cycles from $\kappa_{i}$ and $\kappa_{i}'$,
\textit{cf.} Lemma~\ref{lem:h1bas}:
\begin{gather}
  \label{eq:38}
  \int\nolimits_{\partial S_{j}}\kappa_{H,i} = \delta_{ij}\;,\quad
  \int\nolimits_{\partial S_{j}'}\kappa_{H,i} = 0\;,\quad
  \int\nolimits_{\partial S_{j}'}\kappa_{H,i}' = \delta_{ij}\;,\quad
  \int\nolimits_{\partial S_{j}}\kappa_{H,i}' = 0\;.
\end{gather}

\begin{lemma}
  \label{lem:uniqbas}
  Given interior and exterior Seifert surfaces $S_{i},S_{i}'$, the conditions
  \eqref{eq:38} uniquely determine a symplectic basis $\{\kappa_{H,1},\ldots,\kappa_{H,g},
  \kappa_{H,1}',\ldots,\kappa_{H,g}'\}$ of $\Ch^{1}(\partial\dom)$.
\end{lemma}

\begin{proof}
  If there was another basis complying with \eqref{eq:38}, the differences
  of the basis forms would harmonic 1-forms with vanishing integral over
  \emph{any} cycle. They must vanish identically.
\end{proof}

Given a symplectic basis, we can embark on the canonical construction of
complete  Lagrangian subspaces of $\Ch^{1}(\partial\dom)$ presented in
Remark \ref{rem:ls}. We start
from a partition
\begin{gather}
  \label{eq:39}
  I \cup I' = \{1,\ldots,g\}\quad,\quad I\cap I' = \emptyset\;.
\end{gather}
Owing to Lemma~\ref{lem:4} and \eqref{inter_pairing} the symplectic
pairing $[\cdot,\cdot]_{\partial\dom}$ has the matrix representation
\begin{gather}
  \label{eq:40}
   \begin{bmatrix}
     \mathbf{0}_{g\times{}g}&\bf{I}_{g\times{}g}\\-\bf{I}_{g\times{}g}&\mathbf{0}_{g\times{}g}
   \end{bmatrix}\in\bbR^{2g,2g}\;,
\end{gather}
with respect to the basis
\begin{equation}
  \label{eq:41}
  \Big(\{\kappa_{H,i}\}_{i\in{}I}\cup\{-\kappa_{H,i}'\}_{i\in{}I'}\Big)\cup
  \Big(\{-\kappa_{H,i}'\}_{i\in{}I}\cup\{\kappa_{H,i}\}_{i\in{}I'}\Big)
\end{equation}
of $\Ch^{1}(\partial\dom)$. Thus,
\begin{gather}
  \label{eq:42}
  L_{\Ch} := \operatorname{span}\{\kappa_{H,i}\}_{i\in{}I}\cup\{-\kappa_{H,i}'\}_{i\in{}I'}
\end{gather}
will yield a complete Lagrangian subspace of $\Ch^{1}(\partial\dom)$. By theorems
\ref{theo:closed} and \ref{theo:coclosed}, $L_{\Ch}$ induces self-adjoint
$\curl = \hop \sd$ operators. From Lemma~\ref{Hodge-closed}, Lemma~\ref{Hodge-coclosed}
and \eqref{eq:38} we learn that their domains allow the characterization
\begin{align}
  \label{eq:43}
  \Cd(\curl_{s}) & := \bigl\{\omega\in W^{1}(\sd,\dom):\;
  \sd(i^{\ast}\omega)=0,\;
  \int\nolimits_{\partial S_{j}}\omega = 0,\,j\in I,\;
  \int\nolimits_{\partial S_{j}'}\omega = 0,\,j\in I'\bigr\}\; \\
  \intertext{in the case of closed traces, and}
  \label{eq:43a}
  \Cd(\curl_{s}) & := \bigl\{\omega\in W^{1}(\sd,\dom):\;
  \sd\hop(i^{\ast}\omega)=0,\;
  \int\nolimits_{\partial S_{j}}\omega = 0,\,j\in I,
  \int\nolimits_{\partial S_{j}'}\omega = 0,\,j\in I'\bigr\}\;,
\end{align}
in the case of co-closed traces, respectively. In fact, the choice $I'=\emptyset$
together with closed trace is the one proposed in \cite{YOG90} to obtain a
self-adjoint $\curl$.


\section{Spectral properties}
\label{sec:spectral-properties}

Having constructed self-adjoint versions of the $\curl$ operator, we go on to verify
whether their essential spectrum is confined to $0$ and their eigenfunctions can form
a complete orthonormal system in $\Ltwov$. These are common important features of
self-adjoint partial differential operators.

The following compact embedding result is instrumental in investigating the spectrum
of $\curl_{s}$. Related results can be found in \cite{WEB80} and \cite{PIC84}.

\begin{theorem}[Compact embedding]
  \label{thm:comp}
  The spaces, endowed with the $W^{1}(\sd,\dom)$-norm, 
  \begin{align*}
    X_{0} := & \{\omega\in W^{1}(\sd,\dom):\;\sd^{\ast}\omega = 0,\; i^{\ast}(\hop\omega)=0\}\\
    \text{and}\quad
    X^{\perp} := &
    \{\omega\in W^{1}(\sd,\dom):\;\sd^{\ast}\omega = 0,\; \sd\hop (i^{\ast}\omega)=0\}
  \end{align*}
  are compactly embedded into $L^{2}(\Lambda^{1}(\dom))$.
\end{theorem}

\begin{Remark}
  \label{rem:comp}
  In terms of Euclidean vector proxies these spaces read
  \begin{align*}
    X_{0} = & \{\Vv\in\Hcurl:\; \Div\Vv = 0,\; \gamma_{n}\Vu = 0 \}\;,\\
    X^{\perp} = & 
    \{\Vv\in\Hcurl:\; \Div\Vv = 0,\; \bDiv(\gamma_{t}\Vu) = 0 \}\;\,
  \end{align*}
  where the constraint $\bDiv(\gamma_{t}\Vu) = 0$ should be read as ``orthogonality''
  to $\bgrad H^{\frac{1}{2}}(\partial D)$ in the sense of the Hodge decomposition.
\end{Remark}

\begin{proof}[of Thm.~\ref{thm:comp}]
  The proof will be given for $X^{\perp}$ only. The simpler case of $X_{0}$ draws on
  the same ideas. We are using vector proxy notation, because the proof takes
  us beyond the calculus of differential forms. Note that the inner product
  chosen for the vector proxies does not affect the statement of the theorem.

  A key tool is the so-called regular decomposition theorem that was discovered in
  \cite{BIS87}, consult \cite[Sect.~2.4]{HIP02} for a comprehensive presentation
  including proofs. It asserts that there is $C>0$ depending only on $\dom$
  such that for all $\Vu\in\Hcurl$ there are functions $\Phibf\in(\Hone)^{3}$,
  $\varphi\in\Hone$, with
  \begin{gather}
    \label{eq:56}
    \Vu = \Phibf + \grad\varphi\quad,\quad
    \NHone{\Phibf} + \SNHone{\varphi} \leq C\NHcurl{\Vu}\;.
  \end{gather}
  Let ${(\Vu_{n})}_{n\in\bbN}$ be a bounded sequence in $X^{\perp}$ that is
  \begin{gather}
    \label{eq:57}
    \Div\Vu_n = 0\quad\text{in }\dom \quad\text{and}\quad\bDiv(\gamma_{t}\Vu_n) = 0
    \quad\text{on }\partial\dom\;,\\
    \label{eq:58}
    \exists C>0:\quad \NLtwo{\Vu_{n}} + \NLtwo{\curl\Vu_n} \leq C\;.
  \end{gather}
  Write $\Vu_{n} = \Phibf_{n} + \grad\varphi_{n}$ for the regular decomposition
  according to \eqref{eq:56}. Thus, ${(\Phibf_{n})_{n\in\bbN}}$ is bounded in
  $(\Hone)^{3}$ and, by Rellich's theorem, will possess a sub-sequence that converges
  in $\Ltwov$. We pick the corresponding sub-sequence of $(\Vu_{n})_{n\in\bbN}$
  without changing the notation. 

  Further,
  \begin{align}
    \label{eq:59}
    \Div\Vu_{n} = 0 \quad & \Rightarrow\quad -\Delta\varphi_{n} = \Div\Phibf_{n}\quad
    \text{(bounded in $\Ltwo$)}\;,\\
    \label{eq:60}
    \bDiv(\gamma_{t}\Vu) = 0 \quad  & \Rightarrow\quad
    -\Delta_{\partial\dom}(\gamma\varphi_{n}) = \bDiv (\gamma_{t}\Phibf_{n})\quad
    \text{(bounded in $\Hm[\partial\dom]{-\frac{1}{2}})$}\;.
  \end{align}
  We conclude that $(\gamma\varphi_{n})_{n\in\bbN}$ is bounded in
  $\Hone[\partial\dom]$ and, hence, has a convergente sub-sequence in
  $\Hh[\partial\dom]$ (for which we still use the same notation). The harmonic
  extensions $\widetilde{\varphi}_{n}$ of $\gamma\varphi_{n}$ will converge in
  $\Hone$.

  Finally, the solutions $\widehat{\varphi}_{n}\in\Hone$ of the boundary value problems
  \begin{gather}
    \label{eq:62}
    -\Delta \widehat{\varphi}_{n} = \Div\Phibf_{n}\quad\text{in }\dom\quad,\quad
    \widehat{\varphi}_{n} = 0 \quad\text{on }\partial\dom\;,
  \end{gather}
  will possess a sub-sequence that converges in $\Hone$, as
  $(-\Delta_{\mathrm{Dir}})^{-1}\Ltwo$ is compactly embedded in $\Hone$.  Since
  $\varphi_{n} = \widetilde{\varphi}_{n} + \widehat{\varphi}_{n}$, this provides
  convergence of a subsequence of ${(\Phibf_{n} + \grad\varphi_{n})}_{n\in\bbN}$ in
  $\Ltwov$.
\end{proof}

Let $\curl_{s}:\Cd_{s}\subset L^{2}(\Lambda^{1}(\dom))\mapsto
L^{2}(\Lambda^{1}(\dom))$ be one of the self-adjoint realizations of $\curl$
discussed in the previous section. Recall that we pursued two constructions
based on closed and co-closed traces, respectively.

\begin{Remark} Even if the domain $\Cd_{s}$ of the self-adjoint $\curl_s$ is known
  only up to the contribution of a Lagrangian subspace of $L_\Ch$, we can already
  single out special subspaces of $\Cd_{s}$:
\begin{itemize}
\item [(i)] For the $\curl$ operators based on closed traces, see
  Sect.~\ref{sec:self-adjoint-curl}, in particular Thm.~\ref{theo:closed}, we find
\begin{equation}
  \label{eq:44}
  \sd W^{0}(\sd,\dom) \subset \Cd_{s}\;.
\end{equation}
Indeed, for $\omega\in\sd W^{0}(\sd,\dom)$ there exists $\eta\in{}W^{0}(\sd,\dom)$
with $\omega=\sd\eta$. Due to the trace theorem, $i^*\eta$ belongs to
$W^{-\frac{1}{2}}(\sd,\partial\dom)$.  Consequently, it follows from the commutative
relation \eqref{commutative-relation} that $i^*\omega=\sd{}i^*\eta$ belongs to
$\sd{}W^{-\frac{1}{2}}(\sd,\partial\dom)$. We conclude using \eqref{eq:28}.
\item [(ii)] For the $\curl$ operators based on co-closed traces introduced
  in Sect.~\ref{sec:self-curl-assoc}, it follows that
\begin{equation}
  \label{eq:45}
  \sd W_{0}^{0}(\sd,\dom) \subset \Cd_{s}\;.
\end{equation}
This is immediate from the fact that
\begin{equation}
    \eta\in{}W^{0}_0(\sd,\dom) \hbox{ and }\omega=\sd\eta\hbox{ implies}\quad
    i^*\omega=\sd{}i^*\eta=0\;,
\end{equation}
which means that $\omega$ belongs to $\Cd_s$, see \eqref{eq:30}.
\end{itemize}
\end{Remark}

In the sequel, the kernel of $\curl_{s}$ will be required. We recall
that
\begin{gather*}
  \Kern{\curl_{s}} = \Cd_{s}\cap \Kern{\curl_{\max}}
\end{gather*}
is a closed subspace of $L^{2}(\Lambda^{1}(\dom))$. Moreover, since
$\sd^2=0$ and due to \eqref{eq:44} and \eqref{eq:45}, one has
\begin{gather}
  \label{kernclo}
  \sd W^{0}(\sd,\dom) \subset \Kern{\curl_{s}}\quad\text{in the
  closed
  case},\\
  \label{kerncoclo}
  \sd W_{0}^{0}(\sd,\dom) \subset \Kern{\curl_{s}}\quad\text{in the co-closed case}.
\end{gather}

\begin{lemma}
  \label{lem:bdab}
  The operator $\curl_{s}$ is bounded from below on $\Cd_{s}\cap
  \Kern{\curl_{s}}^{\perp}$:
  \begin{gather*}
    \exists C=C(\dom):\quad
    \N{\omega} \leq C \N{\curl_{s}\omega}
    \quad\forall \omega \in
    \Cd_{s}\cap \Kern{\curl_{s}}^{\perp}\;.
  \end{gather*}
\end{lemma}

\begin{proof}
  The indirect proof will be elaborated for the case of co-closed traces only. The
  same approach will work for closed traces.

  We assume that there is a sequence ${(\omega_{n})}_{n\in\bbN}\subset \Cd_{s}\cap
  \Kern{\curl_{s}}^{\perp} $ such that
  \begin{gather}
    \label{eq:46}
    \N{\omega_{n}} = 1\quad,\quad
    \N{\curl\omega_{n}} \leq n^{-1}\quad\forall n\in\bbN\;.
  \end{gather}
  Since $\omega_{n}\in \Kern{\curl_{s}}^{\perp}$, the inclusion \eqref{kerncoclo}
  implies that $d^{\ast}\omega_{n}=0$. As a consequence of \eqref{eq:46},
  ${(\omega_{n})}_{n\in\bbN}$ is a bounded sequence in $X^{\perp}$.
  Theorem~\ref{thm:comp} tells us that it will possess a subsequence that converges
  in $L^{2}(\Lambda^{1}(\dom))$, again we call it ${(\omega_{n})}_{n\in\bbN}$. Thanks
  to \eqref{eq:46} it will converge in the graph norm on $\Cd_{s}$ and the non-zero
  limit will belong to $\Kern{\curl_{s}}\cap \Kern{\curl_{s}}^{\perp}=\{0\}$. This
  contradicts $\N{\omega_{n}}=1$.
\end{proof}

From Lemma~\ref{lem:bdab} we conclude that the range space
$\Range{\curl_{s}}$ is a closed subspace of
$L^{2}(\Lambda^{1}(\dom))$, which means,
\begin{gather}
  \label{eq:49}
  \Range{\curl_{s}} = \Kern{\curl_{s}}^{\perp}\;.
\end{gather}
Thus, we are lead to consider the symmetric, bijective operator
\begin{gather}
  \label{eq:50}
  \mathsf{C} :=
  \curl_{s}:\Cd_{s}\cap\Kern{\curl_{s}}^{\perp}\subset\Kern{\curl_{s}}^{\perp}
  \mapsto \Kern{\curl_{s}}^{\perp}\;.
\end{gather}
It is an isomorphism, when $\Cd_{s}\cap\Kern{\curl_{s}}^{\perp}$ is equipped with the
graph norm, and $\Kern{\curl_{s}}^{\perp}$ with the $L^{2}(\Lambda^{1}(\dom))$-norm.
Its inverse $\mathsf{C}^{-1}$ is a bounded, self-adjoint operator.
\begin{theorem}
  \label{thm:spec}
  The operator $\curl_{s}$ has a pure point spectrum with $\infty$ as sole accumulation
  point. It possesses a complete $L^{2}$-orthonormal system of eigenfunctions.
\end{theorem}

\begin{proof}
  The inverse operator
  \begin{equation}
  \mathsf{C}^{-1}:\Kern{\curl_{s}}^{\perp}
  \mapsto \Cd_{s}\cap\Kern{\curl_{s}}^{\perp}
  \end{equation}
  is even \emph{compact} as a mapping $L^{2}(\Lambda^{1}(\dom))\mapsto
  L^{2}(\Lambda^{1}(\dom))$. Indeed, due to \eqref{kernclo} and \eqref{kerncoclo} the
  range of $\mathsf{C}^{-1}$ satisfies
\begin{gather}
\Cd_{s}\cap\Kern{\curl_{s}}^{\perp}\subset{}X_0\hbox{ in the closed
case},\\
 \Cd_{s}\cap\Kern{\curl_{s}}^{\perp}\subset{}X^\perp\hbox{ in the
co-closed case}.
\end{gather}
  By Theorem \ref{thm:comp}, the compactness follows.  Riesz-Schauder theory \cite[Sect.~X.5]{YOS80} tells us that,
  except for $0$ its spectrum will be a pure (discrete) point spectrum with zero as accumulation point and it will possess a
  complete orthonormal system of eigenfunctions.

  The formula, see \cite[Thm.~5.10]{WEI80},
  \begin{gather}
    \label{eq:51}
    \lambda^{-1}-\mathsf{C}^{-1} = \lambda^{-1}(\mathsf{C}-\lambda)\mathsf{C}^{-1}
  \end{gather}
  shows that for $\lambda\not=0$,
  \begin{align*}
    \bullet\quad & \lambda^{-1}-\mathsf{C}^{-1}\;\text{bijective}\quad
    \Rightarrow\quad \mathsf{C}-\lambda\;\text{bijective}\;,\\
    \bullet\quad & \Kern{\lambda^{-1}-\mathsf{C}^{-1}} = \Kern{\mathsf{C}-\lambda}\;.
  \end{align*}
  Thus, $\sigma(\mathsf{C}) = (\sigma(\mathsf{C^{-1}})\setminus\{0\})^{-1}$ and the
  eigenfunctions are the same.
\end{proof}

\section{$\curl$ and $\curl\curl$}
\label{sec:curl-curlcurl}

\subsection{Self-adjoint $\curl\curl$ operators}
\label{sec:self-adjo-curlc}

In the context of electromagnetism we mainly encounter the self-adjoint operator
$\curl\curl$. Now we explore its relationship with the $\curl$ operators discussed
before. A metric on $\dom$ and an associated Hodge operator $\hop$ will be taken
for granted.

\begin{definition}
  \label{def:curlcurl}
  A linear operator $\mathsf{S}:\Cd(\mathsf{S})\subset
  L^{2}(\Lambda^{1}(\dom))\mapsto L^{2}(\Lambda^{1}(\dom))$ is a $\curl\curl$
  operator, if and only if $\mathsf{S}$ is a closed extension of the
  operator $\hop\sd\hop\sd$ defined for smooth compactly supported 1-forms.
\end{definition}

Two important extensions of the $\curl\curl$ operator are the maximal and the minimal
extensions:

\begin{lemma}
  \label{lem:curlcurlmin}
  The domain of the \emph{minimal} closed extension $(\curl\curl)_{\min}$ of the
  $\curl\curl$ operator is
  \begin{gather}
    \label{eq:52}
    \Cd_{\min}=\Big\{\omega\in{}W^{1}_0(\sd,D):\;\hop\sd\omega\in{}W^{1}_0(\sd,D)\Big\}\
  \end{gather}
  or, equivalently, in terms of Euclidean vector proxies
  \begin{align*}
     \Cd_{\min}
     = &  \Big\{
     \begin{aligned}[t]
       \Vu\in{}\Ltwov:\;\curl\Vu\in{}\Ltwov,\curl\curl\Vu\in\Ltwov,\\
       \gamma_t(\Vu)=0,\,\text{ and } \gamma_t(\curl{(\Vu)})=0\quad\text{on }\partial\dom
       \Big\}.
     \end{aligned}
   \end{align*}
   The adjoint of $(\curl\curl)_{\min}$ is the \emph{maximal} closed extension
   $(\curl\curl)_{\max}$. It is an extension of the $\curl\curl$ operator with domain
   \begin{equation}
     \Cd_{\max}\;=\;\Cd_1\oplus\Cd_2\;,
   \end{equation}
   with
   \begin{align}
     \label{Domain-curlcurl1}
     \Cd_1= & \Big\{\omega\in{}W^1_0(\sd,\dom):\;\hop\sd\omega \in{}W^1(\sd,\dom)\Big\},\\
     \label{Domain-curlcurl2}
     \Cd_2= & \Big\{\omega\in{}L^{2}(\Lambda^{1}(\dom)):\; \sd\hop\sd\omega=0\Big\}\;.
   \end{align}
\end{lemma}

\begin{proof}
  The domain $\Cd_{\hbox{min}}$ of the minimal closure is straightforward. We recall
  the definition of the domain of the adjoint $\mathsf{T}^{\ast}$ of an
  operator $\mathsf{T}:\Cd(\mathsf{T})\subset H\mapsto H$
  \begin{equation}
    \label{eq:53}
    \Cd( \mathsf{T}^*)\;=\;\Big\{ u\in H:\;\exists C_{u}>0:
    \quad ( u, \mathsf{T} v)_{ H}\;\leqslant\;C_{ u}\;\| v\|_{ H}\quad
    \forall v\in \Cd( \mathsf{T})\Big\}\;.
  \end{equation}
  Let $\Cd_{\max}$ stand for the domain of the adjoint of the minimal $\curl\curl$
  operator.  First we show that
\begin{equation}
  \label{sensesubset}
  \Cd_1\oplus\Cd_2\subset\Cd_{\max}\;.
\end{equation}
Let us consider $\omega\in \Cd_1$ and $\eta\in \Cd_{\min}$. By
integration by parts and the isometry properties of $\hop$ we get
\begin{equation}
  \label{eq:54}
  \int\limits_{\dom}\omega\wedge \sd\hop\sd\eta
  \;=\;\int\limits_{\dom}\sd\hop\sd\omega\wedge\eta \leq
  \N{\sd\hop\sd\omega}\N{\eta}\;.
\end{equation}
This involves $\Cd_1\subset{}\Cd_{\max}$.

Now we consider $\omega\in\Cd_2$. The relation $\sd\hop\sd\omega=0$
has to be understood as
\begin{equation}
  \int\limits_{\dom}\sd\hop\sd\omega\wedge \eta = 0
  \quad\forall\eta\in \Lambda^{1}(\dom)\;\text{smooth, compactly supported}\;.
\end{equation}
As the smooth compactly supported 1-forms are dense in $\Cd_{\min}$
with respect to the topology induced by the norm
\begin{equation}
  \big\|\omega\big\|+\big\|\curl(\omega)\big\|+\big\|\curl(\curl(\omega))\big\|\;,
\end{equation}
it follows that
\begin{equation}
  \int\limits_{D}\omega\wedge\sd\hop\sd\eta =0\quad\forall\eta\in{}\Cd_{\min}\;,
\end{equation}
and, finally, $\Cd_{2}\subset{}\Cd_{\max}$. This confirms \eqref{sensesubset}.

Next, we prove
\begin{equation}
  \label{subsetsense}
    \Cd_{\max}\;\subset\;\Cd_1\oplus \Cd_2\;.
\end{equation}
Pick, $\omega\in \Cd_{\max}$. There exists $\varphi\in L^{2}(\Lambda^{1}(\dom))$ such
that
\begin{equation}
  \label{var-form}
  \int\limits_{D}\omega\wedge \sd\hop\sd\eta;=\;
  \int\limits_{D}\varphi\wedge\hop\eta\quad\forall{}\eta\in\Cd_{\min}\;.
\end{equation}
Since $\sd^{\ast}\varphi=0$ (pick $\eta=\sd \nu$ in
\eqref{var-form}), and $\int_{D}\varphi\wedge\hop\eta_{\Ch}=0$ for
$\eta_{\Ch}\in \Ch^{1}(\dom)$, there exists
$\omega_1\in{}W^1(\sd,\dom)$ satisfying
\begin{equation}
  \left\{
    \begin{array}{lllll}
      \hop\sd\hop\sd\omega_{1} = \varphi & \text{ in } \dom,\\[6pt]
      i^{\ast}\omega_1 = 0&\text{ on } \partial{}D.
    \end{array}
  \right.
\end{equation}
Note that this $\omega_1$ belongs to $\Cd_1$. Then $\omega_2=\omega-\omega_1$
satisfies
\begin{equation}
  \int\limits_{D}(\omega-\omega_1) \wedge
  \sd\hop\sd\eta=0\quad
  \forall \eta \in \Cd_{\min}
  \quad\Longrightarrow{}\quad \sd\hop\sd\omega_2=0\;.
\end{equation}
It follows that $\omega_2\in\Cd_2$. Since $\omega=\omega_1+\omega_2$, we have
proven \eqref{subsetsense}.
\end{proof}

\begin{Remark}
  \label{rem:7}
  The last lemma gives a nice example for
  \begin{displaymath}
    (\mathsf{T}^2)^*\neq{}(\mathsf{T}^2)^*.
  \end{displaymath}
  Indeed, the minimal extension of the formal $\curl\curl$ boils down to
  the squared minimal $\curl$ operator $\curl_{\min}$ with domain $W^1_0(\sd,D)$
  \begin{displaymath}
    (\curl\curl)_{\min}\;=\;\curl_{\min}\;\curl_{\min}
  \end{displaymath}
  The adjoint of $\curl_{\min}$ is the $\curl_{\max}$ operator with domain
  $W^1(\sd,D)$, but
  \begin{displaymath}
    (\curl\curl)_{\max}\;\neq\;\curl_{\max}\;\curl_{\max}\;.
  \end{displaymath}
\end{Remark}

To identify self-adjoint $\curl\curl$ operators we could also rely on the toolkit
of symplectic algebra, using the metric-dependent symplectic pairing
\begin{equation}
  [\omega,\eta]\;=\;
  \int\nolimits_{D}\sd\hop\sd\omega\wedge\eta-\int_{D}\omega\wedge\sd\hop
  \sd\eta\;.
\end{equation}
As before, complete Lagrangian subspaces will give us self-adjoint extensions of
$(\curl\curl)_{\min}$ that are restrictions of $(\curl\curl)_{\max}$. However,
we will not pursue this further.

There are two classical self-adjoint $\curl\curl$ operators that play a central
role in electromagnetic boundary value problems. Their domains are
\begin{align}
  \Cd((\curl\curl)_{\mathrm{Dir}}) = &
  \Big\{\omega\in{}W_{0}^1(\sd,D):\;\hop\sd\omega\in{}W^1(\sd,D)\Big\},\\
  \Cd((\curl\curl)_{\mathrm{Neu}})=
  & \Big\{\omega\in{}W^1(\sd,D):\;\hop\sd\omega\in{}W_0^1(\sd,D)\Big\}\;.
\end{align}
Both can be written as the product of a $\curl$ operator and its adjoint:
\begin{gather}
  \label{eq:55}
  (\curl\curl)_{\mathrm{Dir}} = \curl_{\max}\curl_{\min}\quad,\quad
  (\curl\curl)_{\mathrm{Neu}} = \curl_{\min}\curl_{\max}\;.
\end{gather}
Less familiar self-adjoint $\curl\curl$ operators will emerge from taking the square
of a self-adjoint $\curl$ operator as introduced in Section~\ref{sec:some-examples-self}.

\subsection{Square roots of $\curl\curl$ operators}
\label{sec:square-roots-curlc}

It is natural to ask whether any self-adjoint $\curl\curl$ operator can be obtained
as the square of a self-adjoint $\curl$. We start with reviewing the abstract theory
of square roots of operators, see \cite[Sect.~7.3]{WEI80}.

Let $\mathsf{S}$ be a positive (unbounded) self-adjoint operator on the Hilbert
space $H$. We recall from \cite[Thm.~7.20]{WEI80} that there exists a unique
self-adjoint positive (unbounded) operator $\mathsf{R}$ saytisfying
\begin{equation}
  \label{eq:1sr}
  \mathsf{S}=\mathsf{R}^2, \text{ i.e.}\quad
  \mathsf{D}(\mathsf{S})=\mathsf{D}(\mathsf{R}^2):=\{u\in
  \mathsf{D}(\mathsf{R})\;/\;\mathsf{R}u\in\mathsf{D}(\mathsf{R})\}
  \hbox{ and }\mathsf{S}u=\mathsf{R}^2u\hbox{ if }u\in\mathsf{D}(\mathsf{S})\;.
\end{equation}

\begin{lemma}[domain of square roots]
\label{domain-square-root}
Let $\mathsf{R}_1$ and $\mathsf{R}_2$ be two closed densely defined unbounded
operators on $\sH$ with domains $\mathsf{D}(\mathsf{R}_1)$,
$\mathsf{D}(\mathsf{R}_2)\subset{}\sH$.

If $\mathsf{R}_1^*\;\mathsf{R}_1=\mathsf{R}_2^*\;\mathsf{R}_2$, that is,
\begin{displaymath}
  \mathsf{D}(\mathsf{R}_1^*\;\mathsf{R}_1)\;=\;\mathsf{D}(\mathsf{R}_2^*\;
  \mathsf{R}_2)\quad\hbox{ and }\quad\forall{}u\in\;\mathsf{D}(\mathsf{R}_1^*
  \;\mathsf{R}_1),\;\mathsf{R}_1^*\;\mathsf{R}_1u\;=\;\mathsf{R}_2^*\;\mathsf{R}_2u\;,
\end{displaymath}
then $\mathsf{D}(\mathsf{R}_1)\;=\;\mathsf{D}(\mathsf{R}_2)$.
\end{lemma}

\begin{proof}
  For $i=1,2$, $\mathsf{D}(\mathsf{R}_i)$ equipped with the scalar product
  $(\su,\sv)_i\;=\;(\su,\sv)_\sH+(\mathsf{R}_i\su,\mathsf{R}_i\sv)_{\sH}$ is a
  Hilbert space.\\[6pt]Let us first prove that
  $\mathsf{D}(\mathsf{R}_i^*\mathsf{R}_i)$ is dense in $\mathsf{D}(\mathsf{R}_i)$
  with respect to $(\cdot,\cdot)_i$. We consider
  $\su\in{}\mathsf{D}(\mathsf{R}_i^*\;\mathsf{R}_i)^{\perp}$
\begin{equation}
  \forall{}\sv\in{}\mathsf{D}(\mathsf{R}_i^*\;\mathsf{R}_i),\quad0\;=\;(\su,\sv)_i\;=\;(\su,\sv)_\sH\;+\;(\mathsf{R}_i\su,\mathsf{R}_i\sv)_\sH\;=\;(\su,\sv\;+\;\mathsf{R}_i^*\;\mathsf{R}_iv)_\sH
\end{equation}
As ${\sId}+\mathsf{R}_i^*\;\mathsf{R}_i$ is surjective from
$\mathsf{D}(\mathsf{R}_i^*\;\mathsf{R}_i)$ to $\sH$, see \cite[Theorem 13.31]{RUD73},
$\su$ is equal to zero.

Hence, the spaces $\mathsf{D}(\mathsf{R}_1)$ and $\mathsf{D}(\mathsf{R}_2)$ share the
dense subspace $\mathsf{D}(\mathsf{R}_1^*\;\mathsf{R}_1)=\mathsf{D}(\mathsf{R}_2^*\;\mathsf{R}_2)$.
Moreover, their scalar products coincide on this subset:
\begin{displaymath}
  (\su,\sv)_\sH+(\sR_1\su,\sR_1\sv)_\sH=(\su,\sv\;+\;\sR_1^*\;\sR_1\sv)_\sH=
  (\su,\sv\;+\;\sR_2^*\;\sR_2\sv)_\sH\;=\;(\su,\sv)_\sH\;+\;(\sR_2\su,\sR_2\sv)_\sH.
\end{displaymath}
We conclude using Cauchy sequences.
\end{proof}

Surprisingly, the simple self-adjoint operator $(\curl\curl)_{\mathrm{Dir}}$ does not
have a square root that is a self-adjoint $\curl$:

\begin{lemma}
  \label{lem:nosr}
  The $\curl\curl$ operator $\curl_{\max}\curl_{\min}$
  does not have a square root that is a self-adjoint $\curl$.
\end{lemma}

\begin{proof}
  Let us suppose that $\sT=\curl_{\max}\curl_{\min}$
  admits a $\curl$ self-adjoint square root $\sS$ which implies
  that
  \begin{equation}
    \curl_{\max}\curl_{\max}^*=\curl_{\max}\curl_{\min}= \sT
    =\sS^2 =\sS\;\sS^*\;.
  \end{equation}
  since $\curl_{\max}$ and $\curl_{\min}$ are adjoint and $\sS$ is self-adjoint. Due
  to lemma \ref{domain-square-root}, we have $\sD(\curl_{\max})=\sD(\sS)$ and therefore
  \begin{equation}
    \sS\;=\;\curl_{\max}
  \end{equation}
  since $\sS$ and $\curl_{\max}$ are both
  $\curl$ operators. Clearly, this is not possible since $\curl_{\max}$ is not self-adjoint.
\end{proof}

\begin{Remark}
  \label{rem:Neu}
We remark that the same arguments apply to the operator
$(\curl\curl)_{\mathrm{Neu}}$.
\end{Remark}

\subsection{$\curl\curl\not=\curl\;\curl^*$ is possible}
\label{counterex}

Finally, we would like to show that not all the self-adjoint $\curl\curl$ operators
are of the form $\sR\;\sR^*$ with $\sR$ a $\curl$ operator.

Following an idea of Everitt and Markus ---a similar construction or the Laplacian is introduced in
\cite{EVM05}--- we consider the self-adjoint $\curl\curl$ operator
\begin{equation}
  \sT^0:\Cd(\sT^0)\subset{}\Ltwov\longmapsto \Ltwov,\quad\quad
  \Vu\longmapsto\curl\curl\Vu
\end{equation}
with domain
\begin{equation}
  \Cd(\sT^0)=\Cd_{\min}\oplus\Cd_2\;,
\end{equation}
where $\Cd_{\min}$ and $\Cd_2$ are defined in \eqref{eq:52} and \eqref{Domain-curlcurl2}.

\begin{proposition}
  \label{prop:nocurl}
  There exists no $\curl$ operator $\mathrm{R}$ such
  that
  \begin{equation}
    \label{TTstar}
    \sT^{0}\;=\;\sR\;\sR^*\;.
  \end{equation}
\end{proposition}

\begin{proof}
  Suppose that there exists a $\curl$ operator $\mathsf{R}$
  satisfying \eqref{TTstar}. By definition of the composition of operators one has
  \begin{displaymath}
    \Cd(\sT^{0})\;=\;\Big\{u\in\sD(\sR^*):\; \mathsf{R}u \in\sD(\sR)\Big\}\;.
  \end{displaymath}
  Hence, this implies
  \begin{displaymath}
    \Cd_2\subset\Cd(\sT^{0})\;\subset\;\sD(\sR^*)\;\subset\;W^1(\sd,D)\;.
  \end{displaymath}
  This is not possible since $\Cd_2$ is not a subspace $W^1(\sd,D)$.

  This can be illustrated by means of vector proxies and in the case
  of the unit sphere $D$. Consider the function
  \begin{displaymath}
    \Vu(r,\theta,z)\;=\;\Big(\sum_{n=1}^{+\infty}r^n\;\sin{n\theta}\Big)\;\Ve_z\;,
  \end{displaymath}
  given the cylindrical coordinates. The $\curl$ and $\curl\curl$ of $\Vu$ are
  \begin{gather*}
    \curl\Vu\;=\;\Big(\sum_{n=1}^{+\infty}n\;r^{n-1}\;\cos{n\theta}\Big)\;\Ve_r
    \;-\;\Big(\sum_{n=1}^{+\infty}n\;r^{n-1}\;\sin{n\theta}\Big)\;\Ve_\theta\;,\\
    \curl\curl\,\Vu\;=\;0\;.
  \end{gather*}
  Direct computation leads to
  \begin{displaymath}
    \big\|\Vu\big\|^2\;<\;+\infty\quad\hbox{ and }\quad\big\|\curl\Vu\big\| = +\infty.
  \end{displaymath}
  Hence, this $\Vu$ satisfies $\Vu\in\Cd_2$ but
  $\Vu\notin{}\Hcurl$.
\end{proof}

\begin{Remark}
  \label{rem:9}
    In the same way, we show that there exists
    no $\curl$ operators $\sR_1$ and $\sR_2$ satisfying
\begin{equation}
  \sT^{0}=\sR_1\;\sR_2\;.
\end{equation}
\end{Remark}

\nocite{ARK98,KOT87,GRK04,KOT90a,MAS97,BOT82,CRK02,KOT04,KOT89c,KOT96}

\begin{appendix}
\textbf{Appendix.} Frequently used notations:
\begin{center}%
\begin{longtable}[c]{p{0.2\textwidth}p{0.78\textwidth}}
  $\dom$ & bounded (open) Lipschitz domain in affine space $\bbR^{3}$\\
  $\dom'$ & (compactified) complement $\dom':=\mathbb{R}^{3}\setminus\bar{\dom}$ \\
  $\partial\dom$ & boundary of $\dom$\\
  $\Vn$ & exterior unit normal vector field on $\partial \dom$\\
  $\Vu,\Vv,\ldots$ & vector fields on a three-dimensional domain \\
  $\omega,\eta,\ldots$ & differential forms \\
  $\sv,\su$ & elements of a factor space/trace space of vector  proxies\\
  $\cdot$ & Euclidean inner product in $\bbR^{3}$ \\
  $\times$ & cross product of vectors $\in\bbR^{3}$ \\
  $\mathsf{T},\mathsf{S},\ldots$ & (unbounded) linear operators on a Hilbert space \\
  $\mathsf{T}^{\ast}$ & adjoint operator \\
  $\mathsf{T}_{min}$ & The minimal closure of $\mathsf{T}$ \\
  $\mathsf{T}_{max}$ & The maximal closure of $\mathsf{T}$ \\
  $\Cd(\mathsf{T})$ & domain of definition of the linear operator $\mathsf{T}$\\
  $\Kern{\mathsf{T}}$ & kernel (null space) of linear operator $\mathsf{T}$ \\
  $\Range{\mathsf{T}}$ & range space of an operator $\mathsf{T}$ \\
  $C^{\infty}(\dom)$ & space of infinite differentiable functions on $\dom$ \\
  $\VC^{\infty}(\dom)$ & space of smooth vector fields $(C^{\infty}(\dom))^{3}$\\
  $C^{\infty}_{0}(\dom)$ & functions in $C^{\infty}(\dom)$ with compact support in
  $\dom$\\
  $\VC_{0}^{\infty}(\dom)$ & vector fields in $(C^{\infty}_{0}(\dom))^{3}$\\
  $\Ltwo$ & real Hilbert space of square integrable functions on $\dom$\\
  $\Ltwov$ & square integrable vector fields in $(\Ltwo)^{3}$\\
  $\Hcurl$ & real Hilbert space $\{\Vv\in\Ltwov:\,\curl\Vv\in \Ltwov\}$ with graph
  norm\\
  $\zHcurl$ & closure of $\VC_{0}^{\infty}(\dom)$ in $\Hcurl$ \\
  $\gamma_{t}$ & tangential boundary trace of a vector field \\
  $\gamma_{n}$ & normal component trace of a vector field \\
  $\bgrad$ & surface gradient \\
  $\bscurl$ & scalar valued surface rotation \\
  $\bDiv$ & surface divergence \\
  $\sd$ & exterior derivative of differential forms \\
  $\Lambda^{k}(M)$ & differential $k$-forms on manifold $M$ \\
  $\wedge$ & exterior product of differential forms \\
  $\hop_{g}$ & Hodge operator induced by metric $g$ \\
  $(\cdot,\cdot)_{k,M}$ & inner product on $\Lambda^{k}(M)$ induced
  by a Hodge operator $\hop$ \\
  $L^{2}(\Lambda^{k}(M))$ & Hilbert space of square integrable $k$-forms on $M$ \\
  $\N{\cdot}$ & norm of $L^{2}(\Lambda^{k}(M))$ (``$L^{2}$-norm''):
  $\N{\omega}^{2} := (\omega,\omega)_{k,M}$ \\
  $W^{k}(\sd,\dom)$ & Sobolev space of square integrable $k$-forms with square
  integrable exterior derivative \\
  $W^{k}_{0}(\sd,\dom)$ & completion of compactly supported $k$-forms in
  $W^{k}(\sd,\dom)$ \\
  $i^{\ast}$ & natural trace operator for differential forms \\
  $\left[\cdot,\cdot\right]$ & generic symplectic pairing\\
  $\left[\cdot,\cdot\right]_{M}$ & symplectic pairing of 1-forms on 2-manifold $M$\\
  $\left[\cdot\right]_{\Gamma}$ & jump of trace of a function across 2-manifold $\Gamma$\\
  ${L}^{\symporth}$ & symplectic orthogonal of subspace $L$ of a symplectic space\\
  $\langle\cdot\rangle$ & (relative) homology class of a cycle \\
  $\bbH_i(A;R)$ & $i^{th}$
  homology group of $A$ with coefficients in $R$\\
  $\bbH^i(A;R)$ & $i^{th}$
  co-homology space of $A$ with coefficients in $R$ \\
  $\bbH_i(A,B;R)$ & $i^{th}$ relative homology group of $A$ relative to $B$ with
  coefficients in $R$ \\
  $\bbH^i(A,B;R)$ & relative
  co-homology space of $A$ relative to $B$ with coefficients in $R$ \\
  $\Ch^{1}(\partial\dom)$ & co-homology space of harmonic 1-forms on $\partial\dom$
  \\
  $\mathbb{H}^{1}(\partial\dom)$ & first co-homology factor space of non-exact
  closed 1-forms on $\partial\dom$ \\
  $\mathbb{H}_{1}(\partial\dom)$ & first homology factor space of non-bounding
  cycles on $\partial\dom$ \\
  $\left\langle\cdot\right\rangle$ & selects (relative) homology class of a cycle \\ 
  $\tscurl$ & trace space of $W^{1}(\sd,\dom)$\\ 
  $\mathbf{H}^{-\frac{1}{2}}_{\mathbf{t}}(\bscurl,\partial{}D)$& tangential traces of vector fields
  in $\Hcurl$\\
  $\mathsf{Z}^{-\frac{1}{2}}(\partial\dom)$ & closed 1-forms in $\tscurl$ \\
  $\omega^{0},\omega^{\perp}$ & components of the Hodge decomposition of
  $\omega\in\tscurl$ \\
  \end{longtable}
\end{center}
\end{appendix}


\end{document}